\documentclass[review]{siamart1116}
\nonstopmode

\usepackage[utf8]{inputenc} 
\usepackage{psfrag}
\usepackage{epsfig}
\usepackage{bf1}
\usepackage{amsmath}
\usepackage{amssymb}
\usepackage{amsfonts}
\usepackage{algpseudocode} 
\usepackage{tabularx} 
\usepackage{booktabs} 
\usepackage{multirow}
\usepackage{color}
\usepackage[english]{babel}
\usepackage[labelfont=sc,textfont=it]{caption}
\usepackage{comment}
\usepackage{soul}

\newcommand{\sign}{\textit{sign}}

\usepackage{tikz}
\usetikzlibrary{positioning}

\renewcommand{\vec}[1]{\mathbf{#1}}

\providecommand{\mat}[1]{#1}

\providecommand{\matT}[1]{{#1}^{T}}
\providecommand{\matI}[1]{{#1}^{-1}}

\title{Chronos: A general purpose classical AMG solver for High Performance Computing}
\author{ Giovanni Isotton \footnotemark[1]
    \and Matteo Frigo \footnotemark[1]
    \and Nicol\`o Spiezia \footnotemark[1]
    \and Carlo Janna \footnotemark[1] \footnotemark[2]}

\begin{document}
\maketitle

\renewcommand{\thefootnote}{\fnsymbol{footnote}}
\footnotetext[1] {M$^3$E s.r.l., via Giambellino 7, 35129 Padova, Italy, {\tt e-mail}
g.isotton@m3eweb.it,\,m.frigo@m3eweb.it,\,n.spiezia@m3eweb.it,\,c.janna@m3eweb.it}
\footnotetext[2]{corresponding author}
\renewcommand{\thefootnote}{\arabic{footnote}}

\begin{abstract}
The numerical simulation of the physical systems has become in recent years a
fundamental tool to perform analyses and predictions in several application fields,
spanning from industry to the academy. As far as large scale simulations are concerned,
one of the most computationally expensive task is the solution of linear systems
arising from the discretization of the partial differential equations governing the
physical processes.
This work presents Chronos, a collection of linear algebra functions specifically
designed for the solution of large, sparse linear systems on massively parallel
computers (https://www.m3eweb.it/chronos/). 
Its emphasis is on modern, effective and scalable AMG preconditioners for High Performance
Computing (HPC).
This work describes the numerical algorithms and the main structures of this software suite,
especially from the implementation standpoint. Several numerical results arising from
practical mechanics and fluid dynamics applications with hundreds of millions of unknowns
are addressed and compared with other state-of-the-art linear solvers,
proving Chronos efficiency and robustness.

\end{abstract}


\section{Introduction}
\label{intro}

The solution of linear systems of equations is a central problem in a huge number
of applications in both engineering and science. These problems are particularly
crucial in the simulation of physical processes through the solution of partial
differential equations or system of partial differential equations. In large-scale simulations, 
the solution of linear systems can be the most expensive
task accounting for up to 99\% of the total simulation cost.

In this work, we are interested in developing fast solution algorithms, suitable
for High Performance Computers (HPC), for the linear system:
\begin{equation}
A \vec{x} = \vec{b}
\label{linSys}
\end{equation}
where $A \in \mathbb{R}^{n \times n}$ is the system matrix, $\vec{b}$ and
$\vec{x}\in \mathbb{R}^{n}$ are the right-hand side and solution vector,
respectively, and $n$ is the number of equations.  Although extension
to general matrices is also possible, the present work restricts its focus
to symmetric and positive definite (SPD) matrices which are very common in
most mechanics and fluid dynamics applications.

In current industrial applications, $n$ can easily grow up to few hundreds of million of 
unknowns. On the other side, systems with billions of unknowns have also been
solved in research experiments. The main
difference between industrial problems and these huge research experiments is
that the former are characterized by complex geometries, irregular discretizations
and heterogeneities in the matrix coefficients, while the latter are generally
obtained by successive refinements of regular or quite regular grids.
Despite their smaller size, problems arising from real-world applications
are very challenging and even having large computational resources may not be
enough.

There are several methods to solve~\ref{linSys}, both
direct~\cite{KorGup16,Rouet20,AmeEt19} and
iterative~\cite{SaaVor00,FalSch14,BadMarPri16}, giving excellent performance
on parallel computers. The former are generally preferred in industrial applications
as they are typically more robust and require no experience from the user. The
main downside is that, especially in 3D problems, the matrix factors require a
huge amount of memory thus becoming the limiting factor for large scale simulations.
This work is focused on iterative methods, more specifically on Algebraic Multigrid
(AMG) preconditioning of iterative methods, because these latter present by far less
memory restrictions and are suitable for almost perfectly parallel implementations.
Moreover, in several practical cases, AMG preconditioning guarantees convergence
in a number of iterations that does not depend or only slightly depends on the mesh
size~\cite{TroOosSch01,BreTonBec06,XuZik17}, a property of paramount
importance for the extreme-size simulations that are foreseen in the near future.
The main drawback of AMG preconditioning is that it is still far from being a
black-box method, requiring an experienced user and sometimes a fine
tuning of the set-up parameters. For most AMG solvers, a wrong set-up can easily
lead to slow convergence or overly expensive preconditioners, and, in the worst cases,
even to a failure in the solution~\cite{KorLuGul14}. AMG preconditioners can be
divided into two principal families, classical AMG that are typically more effective
on fluid dynamic problems and aggregation-based AMG which performs better in
solid mechanics.

The Chronos package is a library of iterative methods and AMG preconditioners designed
for high performance platforms to solve severely ill-conditioned problems arising in
real-world industrial applications. To be effective on a wide range of different
applications, Chronos allows for the choice of several options, from the adaptive
generation of the operator near-kernel to the smoother selection, from coarsening
to prolongation, all of this in the framework of the classical AMG method.
In particular, it will be shown that BAMG interpolation makes this AMG extremely
effective also on mechanical problems without the need to use an aggregation based
coarsening.
From the implementation standpoint, Chronos has been developed for HPC, adopting a
distributed sparse matrix storage scheme where smaller CSR blocks are nested into
a global CSR structure. This storage format, together with the adoption of non-blocking
send/receive messages, allows for a high overlap
between communications and computations thus hiding communication latency even
for a relatively small amount of local operations.
Finally, Chronos has a strongly object-oriented design to be readily linked
to other software, to be used as an innermost kernel in more complex approaches, such
as block preconditioners for multiphysics~\cite{BeiPanBen18,FerFraJanCasTch19,FriCasFer19,
RoyJonLemWat20,WatGre20},and to be easily modified to support emerging hardware as GPU and FPGA
\cite{IsoJanBer20,Hagh20,Zach20}.

The algorithms and methods presented in this work are not radically new, but are
rather known algorithms revisited and highly tuned for challenging industrial
problems from various fields. Particular care has been spent in the general
design of the library in order to make it easily maintainable and amenable
of improvements, without sacrificing performance. The benchmarks provided in
the numerical experiments do not derive from the regular discretization of
artificial problems, instead have been collected, also from other
research/industrial groups, with the specific purpose of validating Chronos against the widest
possible selection of test cases.

The remainder of the paper is organized as follows. In the next section, the classical AMG
method will be briefly outlined with a large emphasis on the specific numerical algorithms
implemented to increase effectiveness. In section 3, the design and structure of
the library are accurately described especially from the implementation standpoint.
The performance of Chronos is finally assessed in section 4 on a set of problems
representative of a wide set of real world problems with a comparison with
other state-of-the-art packages. The paper is closed with some concluding remarks
and ideas for future work.

\section{Classical Algebraic Multigrid framework}
\label{sec:1}

In this section, we briefly give an overview of classical AMG and describe,
from a numerical viewpoint, all the AMG components and options, implemented in Chronos.
One of the strengths of this library is that it offers several options for each AMG
component to allow for the user to tune the best combination for any specific
problem.

Any AMG method is generally built on three main components whose interplay gives the
effectiveness of the overall method:
\begin{itemize}
   \item Smoothing, where an inner preconditioner is applied to damp the high-frequency error components;
   \item Coarsening, in which coarse level variables are chosen for the construction of
         the next level;
   \item Interpolation, defining the transfer operator between coarse and fine variables.
\end{itemize}

In Chronos a fourth component, borrowed from the context of bootstrap and adaptive AMG
\cite{BraBraKahLiv11,BreFalMacManMccRug04,BreFalMacManMccRug05,BreFalMacManMccRug06},
is added to the above three and consists in a method to unveil hidden components of
the near kernel of the linear operator whenever they are not a priori available.

As mentioned before, in the present work we are focused on the classical AMG setting,
and below we will briefly recall the basic concepts behind this method, referring
the interested reader to more detailed and rigorous descriptions in the
works~\cite{Stu01,TroOosSch01,XuZik17}.
For the sake of clearness, we restrict this introduction to a two levels only scheme,
as the multilevel version can be readily obtained by recursion.

The first component that has to be set-up in AMG is the smoother, which is a
stationary iterative method responsible for eliminating the error components associated
with large eigenvalues of $\mat{A}$, referred also as the high-frequency
errors. The smoother is generally defined from a rough approximation of
$A^{-1} \simeq M^{-1}$ and its operator is represented by the following equation:
\begin{equation}
\mat{S} = \mat{I} - \omega \matI{M} \mat{A}, 
\label{eq:smoother}
\end{equation}
where $\mat{I}$ is the identity matrix and $\omega$ a relaxation factor to ensure:
\begin{equation}
\omega \rho(\matI{M} \mat{A}) < 2
\end{equation}
see for instance \cite{FraMagMazSpiJan19} for a short explanation.
Generally, the smoother is given by a simple pointwise relaxation method such as (block)
Jacobi or Gauss-Seidel, with the second one often preferred even though its use on
parallel computers is not straightforward. Unlike other AMG packages such as BoomerAMG
\cite{HenYan02} or GAMG \cite{petsc-web-page} where traditional smoothers like
Gauss-Seidel or Chebyshev are selected by default, Chronos implements
the adaptive Factorized Sparse Approximate Inverse (aFSAI)~\cite{JanFerSarGam15} so
that the preconditioning matrix $\matI{M}$ takes the following explicit form:
\begin{equation}
\matI{M} = G^T G
\label{smo}
\end{equation}
with $G$ lower triangular, so that its application simply requires two matrix-vector
products. This choice is dictated by its almost perfect strong scalability and by
its proven robustness in real engineering problems~\cite{BagFraSpiJan17,JanFerGam15}.
Moreover, the cost of aFSAI application is usually much lower than that of Gauss-Seidel
and Chebyshev since its density, ~i.~e. the ratio between the number of non-zeroes of
$\matI{M}$ and $\mat{A}$, is generally $\sim 0.2\div0.4$.

The second component of AMG is the so-called Coarse-Grid Correction (CGC), which is the
$A$-orthogonal projection operation that should take care of the low-frequency components
of the error. To build CGC in classical AMG, the unknowns of a given level are partitioned
into Fine and Coarse (F/C), with those coarse variables becoming the unknowns of the next
level. The choice of coarse variables is a crucial point in the AMG construction, as it
determines both the rate at which the problem size is reduced and the convergence of the
method. Here, we rely on the concept of Strenght of Connection (SoC), i.e., we associate
to each edge of the adjacency graph of $A$ a measure of its relative importance.
Then, using SoC, we rank the graph connections and filter out those deemed less important.
A maximum independent set (MIS) is finally constructed on the filtered SoC graph to determine
coarse variables.

To facilitate explanation, the system matrix is reordered according to this
partitioning of the unknowns with first fine variables and second coarse ones:
\begin{equation}
\mat{A} = 
\begin{bmatrix}
\mat{A}_{ff}  & \mat{A}_{fc} \\
\matT{A}_{fc} & \mat{A}_{cc} \\
\end{bmatrix}
\label{blk_FC}
\end{equation}
with $\mat{A}_{ff}$ and $\mat{A}_{cc}$ square $n_f \times n_f$ and $n_c \times n_c$ matrices,
respectively.
Using this F/C ordering~(\ref{blk_FC}), the prolongation operator $\mat{P}$ is written as:
\begin{equation}
\mat{P} = 
\begin{bmatrix}
\mat{W} \\
\mat{I}
\end{bmatrix},
\label{eq:prolOperator}
\end{equation}
where $\mat{W}$ is a $n_f \times n_c$ matrix containing the weights for coarse-to-fine
variable interpolation. As the system matrix is SPD, the restriction operator $\mat{R}$
is defined through a Galerkin approach as the transpose of $P$, and the coarse level
matrix $\mat{A_c}$ is simply given by the triple matrix product:
\begin{equation}
\mat{A}_c = \matT{P} \mat{A} \mat{P}
\label{eq:AcDef}
\end{equation}
In practice, fast convergence and rapid coarsening, i.e. high F/C ratios,
are always desired, and the construction of effective prolongation operators
is of paramount importance to conciliate these conflicting requirements.

Having defined all the above components, the set-up phase of the two-level
multigrid method is completed and the iteration matrix is given by:
\begin{equation}
\left( \mat{S} \right)^{\nu_{2}}
\left( \mat{I} - \mat{P} \matI{A_c} \matT{P} \mat{A} \right) 
\left( \mat{S} \right)^{\nu_{1}}
\label{eq:twoGridOp}
\end{equation}
with $\nu_{1}$ and $\nu_{2}$ representing the number of smoothing steps
performed before and after the coarse-grid correction, respectively.

\begin{algorithm}[t!]
\caption{\bf AMG Set-up}
\begin{algorithmic}[1]
\Procedure{AMG\_SetUp}{$\mat{A}_k$}
\State Define $\Omega_k$ as the set of the $n_k$ vertices of the adjacency graph of $\mat{A}_k$;
\If {$n_k$ is small enough to allow for a direct factorization}
\State Compute $\mat{A}_k = \mat{L}_k \mat{L}_k^T$;
\Else 
\State Compute $\mat{M}_k$ such that $\mat{M}_k^{-1} \simeq \matI{A}_k$;
\State Define the smoother as $\mat{S}_k = \left(\mat{I}_k -
\omega_k \mat{M}_k^{-1} \mat{A}_k \right)$;
\State Partition $\Omega_k$ into the disjoint sets $\mathcal{C}_k$ and $\mathcal{F}_k$ via coarsening;
\State Compute the prolongation matrix $\mat{P}_k$ from $\mathcal{C}_k$ to $\Omega_k$;
\State Compute the new coarse level matrix $\mat{A}_{k + 1} = \matT{P}_k \mat{A}_k \mat{P}_k$;
\State Call \text{AMG\_SetUp}$\left(\mat{A}_{k+1}\right)$;
\EndIf
\EndProcedure
\end{algorithmic}
\label{algo:cptAMG}
\end{algorithm}

Algorithms~\ref{algo:cptAMG} and~\ref{algo:applyAMGv} briefly report the general AMG set-up
phase and application in a V-cycle, respectively, in a multilevel framework, where it is
conventionally assumed that $A_0 = A$, $\vec{y}_0 = \vec{y}$ and $\vec{z}_0 = \vec{z}$.
Details on all the computational kernels sketched in~\ref{algo:cptAMG} and their parallel
implementation will be discussed in the next sections/subsections.

\begin{algorithm}[t!]
\caption{\bf AMG application in a V-cycle}
\begin{algorithmic}[1]
\Procedure{AMG\_Apply}{$\mat{A}_k$, $\vec{y}_k$, $\vec{z}_k$}
\If{$k$ is the last level} 
\State Solve $\mat{A}_k \vec{z}_k = \vec{y}_k$ using $\mat{L}_k$, the exact Cholesky
factor of $\mat{A}_k$;
\Else 
\State Compute $\vec{s}_k$ by applying $\nu_1$ smoothing steps to
$\mat{A}_k \vec{s}_k = \vec{y}_k$ with $\vec{s}_0 = \vec{0}$;
\State Compute the residual $\vec{r}_k = \vec{y}_k - \mat{A}_k \vec{s}_k$;
\State Restrict the residual to the coarse grid $\vec{r}_{k+1} = \matT{P}_k \vec{r}_k$;
\State Call \text{AMG\_Apply}$\left(\mat{A}_{k+1},\vec{r}_{k + 1},\vec{d}_{k+1} \right)$;
\State Prolongate the correction to the fine grid $\vec{d}_{k} = \mat{P}_k \vec{d}_{k+1}$;
\State Update $\vec{s}_{k} \leftarrow \vec{s}_{k} + \vec{d}_{k}$;
\State Compute $\vec{z}_k$ by applying $\nu_2$ smoothing steps to
       $\mat{A}_k \vec{z}_k = \vec{y}_k$ with $\vec{z}_{0} = \vec{s}_{k}$;
\EndIf
\EndProcedure
\end{algorithmic}
\label{algo:applyAMGv}
\end{algorithm}


\subsection{Unveiling the operator near Kernel}

The kernel (or null space) associated with the homogeneous discretized operator
arising from the most common PDE or systems of PDE is generally a priori known.
For instance it is well-known that the constant vector is the kernel for the
Laplace operator and rigid body modes constitute the kernel for linear
elasticity problems. The information needed to build these spaces, usually
referred to as test spaces in the adaptive AMG terminology, is readily available
to the user from nodal coordinates or other data retrievable from the
discretization. However, the homogeneous operator kernel is only an approximation
of the true near kernel associated with the fully assembled matrix and does not
take into account all the peculiarities of the problem such as boundary conditions
or the strong heterogeneities in the material properties that often arise in real-world problems. 
In many circumstances, a better test space can be obtained by
simply modifying the initial near kernel suggested by the PDE. In the adaptive
AMG literature \cite{BreFalMacManMccRug05,BreFalMacManMccRug06,BraBraKahLiv11,Lee20},
the test space is found by simply running a few smoothing steps over a random test space
or the initial near kernel, whenever available. However, since the near kernel of $A$
is related to the smallest eigenpairs of:
\begin{equation}
A \boldsymbol{\varphi} = \lambda \boldsymbol{\varphi}
\label{eigP}
\end{equation}
a better way to extract an effective test space could be by relying on an iterative
eigensolver. In the present implementation, we opt for the simultaneous Rayleigh quotient
minimization (SRQM) \cite{BerMarPin06,FraMagMazSpiJan19} whose cost per iteration is only
slightly higher than a smoothing step. By contrast, SRQM can provide a much better
approximation of the smallest eigenpairs especially if a good preconditioner is provided.
Since an approximation of $A^{-1}$ is already available through the smoother, we simply
reuse the previously computed $M^{-1}$ inside the SRQM iteration.

From a theoretical standpoint, instead of solving (\ref{eigP}), the test space
should be computed by solving the generalized eigenproblem:
\begin{equation}
A \boldsymbol{\varphi} = \lambda M \boldsymbol{\varphi}
\label{GeigP}
\end{equation}
However, the SRQM solution to (\ref{GeigP}) needs the multiplication of $M$ by a
vector which, due to our choice of $M$ (\ref{smo}), would result in a forward and
backward triangular solve whose parallelization may represent an algorithmic
bottleneck.

Unfortunately, extracting with high accuracy the eigenpairs of (\ref{eigP}) is 
generally more expensive than solving the original linear system (\ref{linSys}). 
For this reason, to limit the set-up cost, we only approximately solve (\ref{eigP})
with a predetermined and small number of SRQM iterations. This simple strategy usually
gives satisfactory results, whenever an initial test space is not available or
boundary conditions and heterogeneity exert a strong influence, such as in geomechanical
problems. Another appealing idea, though not explored in this work, is bootstrapping
\cite{BraBraKahLiv11,BreFalMacManMccRug05,BreFalMacManMccRug06}, which consists
in computing a relatively cheap AMG preconditioner from a tentative test space,
and then using AMG itself to better uncover the near null space and rebuild a more
effective AMG. It will be shown in section~\ref{sec:3} how the object-oriented
implementation of Chronos allows for easily using such an approach through
simple calls to high-level functions.

Operatively, once the test space is found, we compute an orthonormal basis of it
and collect the basis vectors into a (skinny) matrix $V$ that may be subsequently
used eventually for the strength of connection and the prolongation.

\subsection{Strength of Connection}

The construction of the coarse problem in \\ Chronos is based on the definition of a
SoC matrix, that is used to filter-out weak connections from the
adjacency graph of $A$. There are three different SoC definitions available through
the library:
\begin{enumerate}
   \item Classical strength of connection:
      \begin{equation}
      s_{ij} = \frac{-a_{ij}}{\max (\min_{j \ne i} a_{ij}, \min_{j \ne i} a_{ji})}
      \label{cla}
      \end{equation}
   \item Strength of connection based on strong couplings:
      \begin{equation}
      s_{ij} = \frac{|a_{ij}|}{\sqrt{a_{ii}a_{jj}}}
      \label{sco}
      \end{equation}
   \item Affinity-based strength of connection:
      \begin{equation}
      s_{ij} = \frac{(\sum_k v_{ik} v_{jk})^2}{(\sum_k v_{ik}^2) (\sum_k v_{jk}^2)}
      \label{aff}
      \end{equation}
\end{enumerate}
where $s_{ij}$ denotes the SoC between node $i$ and $j$ and $a_{ij}$ and $v_{ij}$ denote
the entries in row $i$ and column $j$ of the matrices $A$ and $V$, respectively.
SoC (\ref{cla}) is particularly effective for Poisson-like problems where the system
matrix is close to an M-matrix. SoC (\ref{sco}) is generally used in smoothed aggregation
AMG \cite{Van96} and usually gives good results in structural problems. Finally,
SoC (\ref{aff}) has been introduced in \cite{LivBra12} and, though requiring a rather
expensive computation, it is able to accurately capture anisotropies as is shown in
\cite{PalFraJan19}.

After SoC is computed for every pair of nodes, weak connections are eliminated
to determine a Maximum Independent Set (MIS) of nodes that will become coarse nodes 
in the next level. The more aggressively the connections are eliminated, the higher
number of nodes are left in the next level. There are two ways of controlling SoC
filtering in Chronos:
\begin{enumerate}
   \item by a threshold, the traditional way of filtering, where we simply
   drop connections with strength below a given threshold $\theta$;
   \item prescribing an average number of connections per node.
\end{enumerate}
On one side, guaranteeing an average number of connections per node is trickier, 
since it requires a preliminary sorting of all the SoC. On the other, it ensures a
more regular grid coarsening through levels with an almost constant coarsening
ratio. Moreover, in affinity-based SoC, the strength values usually lie in a narrow
interval close to unity so that a proper choice of the drop threshold is almost
impossible.

Finally, MIS construction is performed by using the PMIS strategy introduced
in \cite{DeSYanHey06} which is a perfectly parallel algorithm giving generally
rise to lower complexities than the classical Ruge-St\"{u}ben coarsening
\cite{DeSFalNolYan08}. Using this more aggressive coarsening method requires some
special care in the interpolation construction, as we will see in the next section.

\subsection{Interpolation}
\label{subsec:Int}

Providing a good interpolation operator is crucial for an effective AMG method.
We recall that the prolongation operator $P$ should satisfy:
\begin{equation}
\mathcal{V} \subseteq \mbox{range}(P)
\end{equation}
where $\mathcal{V}$ is the near-kernel of $A$ or, more precisely, for a coarse space
of given size $n_c$ the optimal two-level prolongation as stated in
\cite{XuZik17,BraCaoKahFalHu18} should be such that:
\begin{equation}
\mbox{span}(\vec{v}_i) = \mbox{range}(P)
\end{equation}
where $\vec{v}_i$ are the eigenvectors associated with the smallest $n_c$ eigenvalues
of the generalized eigenproblem (\ref{GeigP}). To this aim, depending on the problem,
we use two different strategies.

If a test space is available or it is relatively cheap
to obtain a reasonable approximation of the near-kernel, then the so-called BAMG
approach is used \cite{BraBraKahLiv11}, where the weights of prolongation $w_{ij}$,
i.e., the entries of the $W$ block in~(\ref{eq:prolOperator}), are found through a least square
minimization:
\begin{equation}
w_{ij} = \underset{j \in C_i}{\mbox{argmin}} \;
         \| \vec{v}_i - \sum_{j \in C_i} w_{ij} \vec{v}_j\|^2 \qquad i = 1,\dots,n
\label{minBAMG}
\end{equation}
with $\vec{v}_k$ the $k$-th row of $V$ and $C_i$ the interpolatory set for $i$. In practice,
it is our experience that to have an effective prolongation, the norm
$\| \vec{v}_i - \sum_{j \in C_i} w_{ij} \vec{v}_j\|$ must be reduced to zero and,
in the general case, this can be accomplished only if the cardinality of $C_i$, $|C_i|$,
is equal or larger than $n_t$, the number of test vectors. To guarantee an exact
interpolation, it is often necessary to use neighbors at a distance larger than
one, especially when dealing with systems of PDEs, with a consequent increase
of the overall operator complexity.
Moreover, it may happen that, even if $|C_i| \ge n_t$, some of the vectors
$\vec{v}_k$ are almost parallel and high conditioning of $\varPhi$, the dense matrix
defined below, may produce large jumps in the weights. In turn, large jumps in $P$
introduces high frequencies in the next level operator that the smoother hardly handles.
To overcome these difficulties, we adopt an adaptive procedure to compute our BAMG
interpolation. More in detail, let us define $\varPhi$ the matrix whose entries
$\varphi_{ij}$ correspond to the $j$-th component of the $i$-th test vector $\vec{v}_i$,
for any $j$ in the interpolatory set.
For each fine node $i \in \mathcal{F}$ to be interpolated, we start by including
in the interpolatory set all its coarse neighbors at a distance no larger than
$l_{\mbox{min}}$, and select a proper basis for $\varPhi$ by using a maxvol
algorithm \cite{Knu85,Gor08}. If either the relative residual:
\begin{equation}
r_i = \frac{\| \vec{\varphi}_i - \varPhi \vec{w}_i  \|}{\| \vec{\varphi}_i \|}
\end{equation}
or the norm of the weights, $\| \vec{w}_i \|$, are larger than the user-defined
thresholds $\epsilon$ and $\mu$, respectively, then we extend by one the
interpolation distance. We keep on increasing the interpolation distance up to
$l_{\mbox{max}}$ to limit the computational cost.  
This procedure, which is briefly sketched in Algorithm \ref{algo:BAMG}, 
though slightly expensive, allows to compute an accurate and smooth
prolongation without impacting too much on the operator complexity. In fact, including
in $C_i$ all the coarse nodes within a priori selected interpolation distance
usually leads to a more complex operator since several fine nodes may be interpolated
with excessively large support. Moreover, limiting the number of non-zeroes in the
rows of $P$ has the additional advantage that it is possible to perform prolongation
smoothing without an exponential growth of the operator complexity. Prolongation 
smoothing is a very common practice in solving elasticity problems with
aggregation-based AMG and numerical results will show how it can beneficial also
in the context of classical AMG.

\begin{algorithm}[t!]
\caption{\bf BAMG prolongation adaptive set-up}
\begin{algorithmic}[1]
\Procedure{BAMG\_Prolongation}{$\mat{S}$, $V$, $l_{\mbox{min}}$, $l_{\mbox{max}}$, $\epsilon$, $\mu$}
\ForAll {$i \in \mathcal{C}$};
\State Set $l = l_{\mbox{min}}$;
\State Set $\vec{r}_i = \vec{\varphi}_i$;
\While {$l \le l_{\mbox{max}}$ {\bf and} ( $r_i > \epsilon$ {\bf or} $\vec{w}_i > \mu$ ) }
\State Include in $C_i$ all the coarse nodes at a distance at most $l$;
\State Collect all the $\varphi_k$ such that $k \in C_i$;
\State Select from $\varphi_k$ a {\em maxvol} basis $\varPhi_i$;
\State Find the vector of weights $\vec{w}_i$ by minimizing $\| \varPhi_i \vec{w}_i - \vec{\varphi}_i\|$;
\State $l = l + 1$;
\EndWhile
\EndFor
\EndProcedure
\end{algorithmic}
\label{algo:BAMG}
\end{algorithm}
%

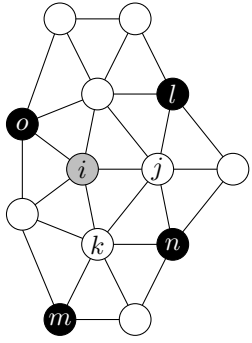
\begin{figure}
\begin{center}
\begin{tikzpicture}

    \draw (0.5,0,0) -- (1.5,0,0); 
    \draw (0.5,0,0) -- (0,1.2,0); 
    \draw (1.5,0,0) -- (1.0,1.0,0); 
    \draw (0,1.4,0) -- (1.0,1.0,0); 
    \draw (0.5,0,0) -- (1.0,1.0,0); 
    \draw (1.0,1.0,0) -- (2.0,1.0,0); 
    \draw (1.5,0.0,0) -- (2.0,1.0,0); 

    \draw (0,1.4,0) -- (0.8,2.0,0); 
    \draw (0,1.4,0) -- (0,2.6,0); 
    \draw (1.0,1.0,0) -- (0.8,2.0,0); 
    \draw (1.0,1.0,0) -- (1.8,2.0,0); 
    \draw (2.0,1.0,0) -- (1.8,2.0,0); 
    \draw (2.0,1.0,0) -- (2.8,2.0,0); 
    \draw (1.8,2.0,0) -- (2.8,2.0,0); 
    \draw (1.8,2.0,0) -- (0.8,2.0,0); 

    \draw (0.8,2.0,0) -- (0.0,2.6,0); 
    \draw (0.8,2.0,0) -- (1.0,3.0,0); 
    \draw (1.8,2.0,0) -- (1.0,3.0,0); 
    \draw (1.8,2.0,0) -- (2.0,3.0,0); 
    \draw (2.8,2.0,0) -- (2.0,3.0,0); 
    \draw (1.0,3.0,0) -- (2.0,3.0,0); 
    \draw (0.0,2.6,0) -- (1.0,3.0,0); 

    \draw (0.0,2.6,0) -- (0.5,4.0,0); 
    \draw (1.0,3.0,0) -- (0.5,4.0,0); 
    \draw (1.0,3.0,0) -- (1.5,4.0,0); 
    \draw (2.0,3.0,0) -- (1.5,4.0,0); 
    \draw (0.5,4.0,0) -- (1.5,4.0,0); 

    \draw[fill=black] (0.5,0,0) circle (0.6em);
    \draw[fill=white] (1.5,0,0) circle (0.6em) ;
    \draw[fill=white] (1.0,1.0,0) circle (0.6em);
    \draw[fill=black] (2.0,1.0,0) circle (0.6em);
    \draw[fill=white] (0,1.4,0) circle (0.6em);
    \draw[fill=lightgray] (0.8,2.0,0) circle (0.6em);
    \draw[fill=white] (1.8,2.0,0) circle (0.6em);
    \draw[fill=white] (2.8,2.0,0) circle (0.6em);
    \draw[fill=black] (0.0,2.6,0) circle (0.6em);
    \draw[fill=white] (1.0,3.0,0) circle (0.6em);
    \draw[fill=black] (2.0,3.0,0) circle (0.6em);
    \draw[fill=white] (0.5,4.0,0) circle (0.6em);
    \draw[fill=white] (1.5,4.0,0) circle (0.6em);

    \node at (0.8,2.0,0.0) {$i$};
    \node at (1.8,2.0,0.0) {$j$};
    \node at (1.0,1.0,0.0) {$k$};

    \node at (0.5,0.0,0.0) {\textcolor{white}{$m$}};
    \node at (2.0,1.0,0.0) {\textcolor{white}{$n$}};
    \node at (2.0,3.0,0.0) {\textcolor{white}{$l$}};
    \node at (0,2.6,0.0) {\textcolor{white}{$o$}};
\end{tikzpicture}
\end{center}
\caption{Example of the interpolatory points. The gray point is the point to be
interpolated, black points are C-points and white points are F-points.}
\label{fig:intset}
\end{figure}

On the other hand, when there is no explicit knowledge of the test space or when the matrix
at hand arises from the discretization of a Poisson-like problem, Chronos can also
rely on more classical interpolation schemes. Below we briefly recall the expressions
of some well-known interpolation formulas. First, using the concept of strength of
connection, we define the following sets:
\begin{itemize}
   \item $N_i = \{ j \; | \; a_{ij} \ne 0 \}$, the set of direct neighbours of $i$;
   \item $S_i = \{ j \in N_i \; | \; j \text{ strongly influences } i \}$, the set of
         strongly connected neighbours of $i$;
   \item $F_i^S = F \cap S_i$, the set of strongly connected fine neighbors of $i$;
   \item $C_i^S = C \cap S_i$, the set of strongly connected coarse neighbors of $i$;
   \item $N_i^W = N_i \setminus ( F_i^S \cup C_i^S )$, the set of weackly connected 
         neighbors of $i$.
\end{itemize}
A generally accurate distance-one interpolation formula, introduced in~\cite{RugStu87},
is the classical interpolation. 
Unlike other distance-one formulas, here, the interpolation takes care of the contribution 
from strongly influencing points $F_i^S$, and the expression for the interpolation weight is given by:

\begin{equation}
w_{ij} = - \frac{1}{a_{ii}+\sum_{k \in N_i^w \cup F_i^{S*}} a_{ik}} \left ( \sum_{k \in F_i^s \setminus F_i^{S*}} 
           \frac{a_{ik}\bar{a}_{kj}}{\sum_{m \in C_i^s} \bar{a}_{km}}  \right ), \qquad  j\in C_i^S,
\end{equation}
where:
\begin{equation}
\bar{a}_{ij}=
\begin{cases}
0  & \text{ if } \sign(a_{ij}) = \sign(a_{ii}) \\
a_{ij} & \text{ otherwise }
\end{cases}
\end{equation}
It is worth noting that the original formula proposed in~\cite{RugStu87} is here
corrected accordingly with the modification introduced in~\cite{HenYan02}
where the set of strongly connected neighbors $F_i^{S*}$, that are F-points but
do not have a common C-point, are subtracted to the fine strong neighbors $F_i^S$.
This modification of the interpolation formula is needed to avoid that the term
$\sum_{m \in C_i^s} \bar{a}_{km}$ vanishes.
Indeed, using the PMIS-coarsening method no longer guarantees that two strongly
connected F-points are interpolated by a common C-point.
However, even if for a large class of problems the classical interpolation is very effective, 
it can lose efficiency for challenging problems such as rotated anisotropies or problems
with large discontinuities.
Indeed, in these cases, the convergence of the AMG accelerated by a Krylov subspace method
can deteriorate, losing scalability and effectiveness.   
Hence, some more advanced interpolation formulas are required to overcome these difficulties.
Particularly, we have to adopt some long-range interpolation strategies.
A widely used interpolation strategy, mainly for very challenging Poisson-like problems,
is the Extended+i interpolation.
This interpolation formula is obtained extending the interpolatory set including C-points
that are distance two away from the F-point we are considering. 
Furthermore, not only the connections from strong fine neighbors to points of the interpolatory set are considered, but also the connections from the 
fine neighbors to the fine point to be interpolated itself.
Hence, denoting with $\hat{C}_i=C_i \cup \bigcup_{j \in F_i^S} C_j$ the set of distance-two coarse nodes,
the interpolation Extended plus i formula takes the following form:

\begin{equation}
w_{ij} = - \frac{1}{\tilde{a}_{ii}} \left ( a_{ij} + \sum_{k \in F_i^s} 
           \frac{a_{ik} \bar{a}_{kj}} { \sum_{l \in \hat{C}_{i\cup \{i\}}} \bar{a}_{kl}}  \right ), \qquad  j\in \hat{C}_i
\end{equation}
with
\begin{equation}
\tilde{a}_{ii} = a_{ii} + \sum_{n \in N_i^w \setminus \hat{C}_i} a_{in} + \sum_{k \in F_i^s} a_{ik} 
\frac{\bar{a}_{ki}}{\sum_{l \in \hat{C}_i \cup {i}} \bar{a}_{kl}}.
\end{equation}
This extended+i interpolation remedies many problems that occur with the distance-one
interpolation formula and it provides better weight coefficients compared with other
distance-two interpolation formulas.
However, unlike the distance-one methods, they lead to much larger
operator complexities. A possible way to reduce the complexities
without or mildly affecting the convergence rate of the iterative scheme is to consider
a different interpolatory set, i.e., an interpolatory set larger than the distance-one
set $C_i^S$, but smaller than the distance-two $\hat{C}_i$.
The idea is to consider an interpolatory set that only extends $C_i^S$ for strong
F-F connections without a common C-point, since in the other cases the point $i$ is
already surrounded by interpolatory points belonging to $C_i^S$.
A crucial point is how to extend $C_i^S$, ensuring as much as possible the quality of
the interpolation operator. 
Here, we propose to enrich the set $C_i^S$ taking into account the minimum number of
distance-two coarse nodes such as to guarantee that each F-F strong connection has at
least a common C-point.  
To better explain this idea, let us consider the example in Figure~\ref{fig:intset}. 

Notice that using the classical interpolation, the interpolatory set would be $C_i^S=\{o\}$ 
and there would be two fine neighbors of $i$, $j$ and $k$, that do not share a
C-point with $i$.  
On the other hand, using the extended plus i interpolation, we would have that
$\hat{C}_i=\{m, n, l, o \}$ and each F-node, strongly connected with $i$, would share
at least one $C$-node of the interpolatory set $\hat{C}_i$.
However, to guarantee this last condition, it would be sufficient to only include the node $n$
to the set $C_i^S$, so that the extended interpolatory set would become $\hat{C}^h_i=\{o, n\}$.
It is worth noting that the new points included in the set $\hat{C}^h_i$ are the minimum
number of coarse point necessary to guarantee that each F-node strongly connected with
the point to be interpolated, is also strongly connected with a C-node belonging to the
interpolatory set.
In other words, we extend the set $C_i^S$ by including the maximum independent set of
distance-two C-nodes such that accomplish the above condition.
Algorithms~\ref{alg:extintset} gives a general description of the procedure used to
make the extended interpolatory set. 

\begin{algorithm}[t!]
\caption{Computation of extended interpolation set $\hat{C}^h_i$}
\begin{algorithmic}[1]
\State{Set the initial interpolatory set $\hat{C}^h_i=C^S_i$}
\State{Compute the initial set $F^{'}$, such that:\\ 
$\qquad F^{'}=\{j \in F^S_i \;  | \; j \text{ is not strongly connected with at least one node in } C_i^S \}$}
\State{Compute the initial set of $C^{''}$, i.e., the set of the distance-two coarse nodes strongly connected with a $F^{'}$-point}
\State{Compute the vertex degree of each element in $C^{''}$ by taking into account only the connections with $F^{'}$}
\While{$F^{'}\ne \emptyset$}
\State{Choose the node with maximum degree in $C^{''}$ and add it to $\hat{C}^h_i$ }
\State{Update the set $F^{'}$}
\State{Update the set $C^{''}$}
\EndWhile
\end{algorithmic}
\label{alg:extintset}
\end{algorithm}


\subsection{Filtering}

One problem that may affect AMG methods, especially in parallel implementation,
is the excessive stencil growth occurring in lower levels. This drawback is even
more pronounced if long-range interpolation or prolongation smoothing is used.
Some authors have explored interesting solutions to reduce AMG complexity without
detrimental effects on convergence \cite{FalSch14,BieFalGroOlsSch16}. 
Simply eliminating small entries from the operators, as is done for instance
with ILU or some approximate inverse preconditioners, may completely harm the
effectiveness of CGC. This happens because removal of small entries from
$P_k$ or $A_{k+1} = P_k^T A_k P_k$ may induce a representation of the near
kernel of $A$ which is not accurate enough for AMG.

To overcome this problem, the authors in \cite{FalSch14} propose to
compensate the action of eliminated entries through a sort of stencil collapsing
to guarantee that the filtered operator, say $\tilde{A}_{k+1}$, exerts on the near
kernel the same action of $A$:
\begin{equation}
\tilde{A}_{k+1} W = A_{k+1} W
\label{action}
\end{equation}
with $W$ a matrix representation of the near kernel. While it is relatively simple
to enforce condition (\ref{action}) for one dimensional near kernels, it is not
immediate to accommodate the action on several vectors at the same time.
With multiple vectors, first the smallest entries of $A_{k+1}$ are dropped to determine
the pattern of $\tilde{A}_{k+1}$, then a correction to $\tilde{A}_{k+1}$, $\Delta_{k+1}$,
is computed by using least squares on:
\begin{equation}
\| (A_{k+1} - \tilde{A}_{k+1}) W = \Delta_{k+1} W \|_2
\end{equation}
More in detail, $\tilde{A}_{k+1}$ is computed row-wisely such that the absolute norm
of each row is a given percentage $\rho$ of the norm of the original one and then
the compensation is computed for the same row. We use the same procedure
on the prolongation operator $P_k$ with the only exception that instead of $W$ we use
its injection in the coarse space. Operatively, the test space $V$ is used when 
available while in cases where $V$ is not computed, such as in Poisson problems,
we simply replace $V$ with a unitary vector.

Finally, we observe that $\tilde{A}_{k+1}$ is no more guaranteed to be SPD and,
especially when an aggressive dropping is enforced, the use of a non-symmetric
Krylov solver, such as GMRES \cite{SaaSch86} or BiCGstab \cite{Vor92}, is often needed
instead of CG. Obviously, such care is not needed when only the prolongation is filtered
as $\tilde{P}_k^T A \tilde{P}_k$ is always SPD for any choice of $\rho$.

\section{Library description}
\label{sec:3}

The Chronos software package is a collection of classes and functions that implements
linear algebra algorithms for distributed memory parallel computers.
The library is written in C++, and Message Passing Interface (MPI) and OpenMP directives
were used for communication among processes and multithread execution, respectively.
The hybrid MPI-OpenMP implementation is more flexible in the use of modern
computing resources and it is generally more efficient than pure MPI due to its
better exploitation of fine-grained parallelism.

Chronos has been developed using the potential of Object-Oriented Programming (OOP).
The abstraction introduced through the OOP allows for using the same distributed
matrix object to represent a linear system, a smoother, an AMG hierarchy
or a preconditioner itself.
Another advantage of this developed approach is the possibility to use simpler classes
to derive more advanced elements, as block preconditioners.
Moreover, whatever the type of preconditioner, the same iterative methods can be used
for the linear system or eigenproblem solution.

In addition, the modular structure allows to easily integrate the CPU kernels
with Graphics Processing Units (GPU) and Field Programmable Gate Array (FPGA)
kernels leaving the overall architecture of the library unchanged, making
Chronos a potentially multi-platform software.
A hybrid CPU-GPU version is already under development and preliminary performances
are encouraging~\cite{IsoJanBer20}.

A brief description of the main classes is reported in the next subsections.

\subsection{Main classes}

\begin{figure*}[h]
\centering
\includegraphics[width=1.0\textwidth]{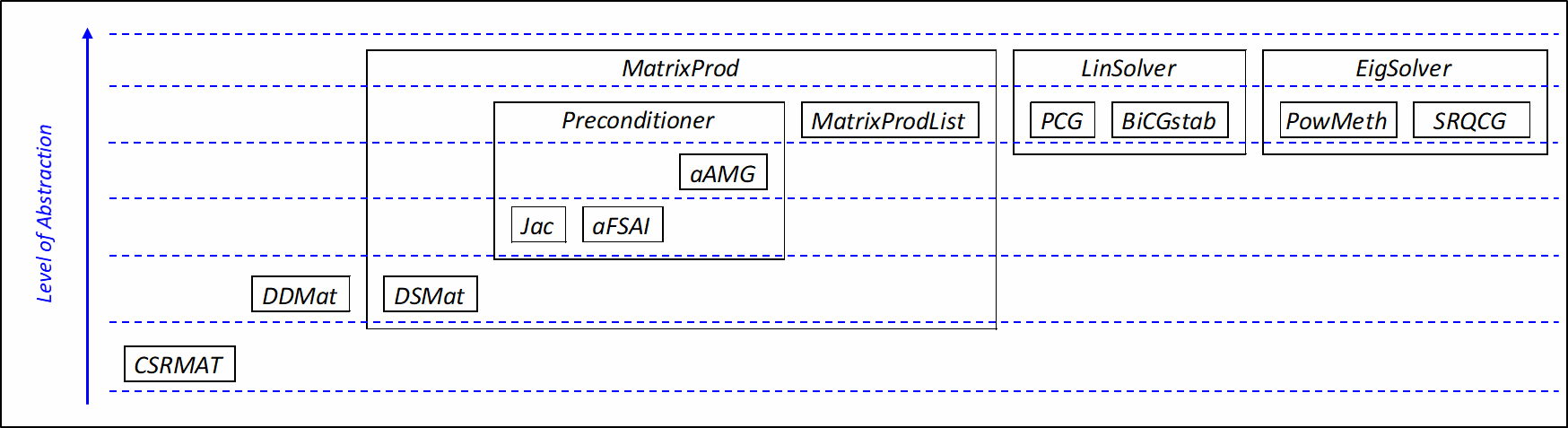}
\caption{Chronos main classes and hierarchies.\label{fig-classes}}
\end{figure*}

The level of abstraction and the hierarchy of the main classes are sketched
in Figure~\ref{fig-classes}.
All these classes are exposed to the user to access the full range of Chronos
functionalities.

The Distributed Dense Matrices (DDMat) and Distributed Sparse Matrices (DSMat) are managed by
the $DDMat$ and the $DSMat$ classes, respectively.
Both DDMat and DSMat storage schemes require the matrix to be subdivided into $n_p$
horizontal stripes of consecutive rows, where $n_p$ is the number of active MPI processes.
In the DDMat, each stripe is stored row-wisely among the process to guarantee
better access in memory during multiplication operations.
This makes the DDMat very efficient for linear systems with multiple
right-hand-sides and eigenproblems, and distributed vectors are stored as
one-column DDMat.
In the DSMat each stripe is subdivided into an array of Compact Sparse Row (CSR)
matrices.  The CSR format is the Chronos standard format for shared sparse matrices
and their management is demanded to $CSRMAT$ class.
The DSMat storage scheme adopted in Chronos is very effective in both the
preconditioner computation and the SpMV product because it allows a large
superposition between communication and computation, and it is described in detail
in the next subsection.

The $Preconditioner$ class manages the approximation of the inverse
of a Distributed Sparse Matrix at the highest level of abstraction.
It requires in input a Distributed Sparse Matrix as
$DSMat$-type object and an optional test space as a $DDMat$-type object.
The classes derived from $Preconditioner$ are $Jac$, $aFSAI$ and $aAMG$ that manage
the preconditioners of Jacobian-type, adaptive-FSAI-type, and aAMG-type, respectively.
A useful feature is that each of these classes can be used as a smoother in the AMG.

Both the $DSMat$ and $Preconditioner$ classes are derived from the $MatrixProd$ class which
manages the Sparse-Matrix-by-Vector product (SpMV) at the highest level of abstraction.
The SpMV is the most expensive operation in any preconditioned iterative solver
and its management has defined the design of the whole library.
With reference to Figure~\ref{fig-classes}, the $MatrixProd$ class leads
the Chronos structure together with the iterative solvers.
Furthermore, more general $MatrixProd$ elements can be readily built using the
$MatrixProdList$ class, that manages an implicit $MatrixProd$ object defined as a product
of a sequence of $MatrixProd$ objects ordered into a list.

At the top of the hierarchy pyramid, there are also the solvers for linear systems and
eigenproblem, $LinSolver$ and $EigSolver$, respectively.
The $LinSolver$ manages the Krylov methods for linear system solution, it requires in
an input preconditioner and a linear system as $MatrixProd$-type objects,
a right-hand-side as a $DDMat$-type object and an optional initial solution as a
$DDMat$-type object.
The classes derived from $LinSolver$ are currently $PCG$ and $BiCGstab$ that manage
the Preconditioned Conjugate Gradient (PCG) iterative method
and the Preconditioned Biconjugate Gradient Stabilized (BiCGstab) iterative method,
respectively.

Finally, the $EigSolver$ manages the Krylov methods for the eigenproblem solution, it
requires an optional input preconditioner and a linear system as $MatrixProd$-type
objects and an initial eigenspace as a $DDMat$-type object.
The two classes derived from $EigSolver$ are currently $PowMeth$ and $SRQCG$, implementing
the Power Method and the Simultaneous Rayleigh Quotient Minimization iterative methods,
respectively.

\subsection{Distributed Sparse Matrix Storage Scheme}
\begin{figure*}[h]
\centering
\includegraphics[width=0.5\textwidth]{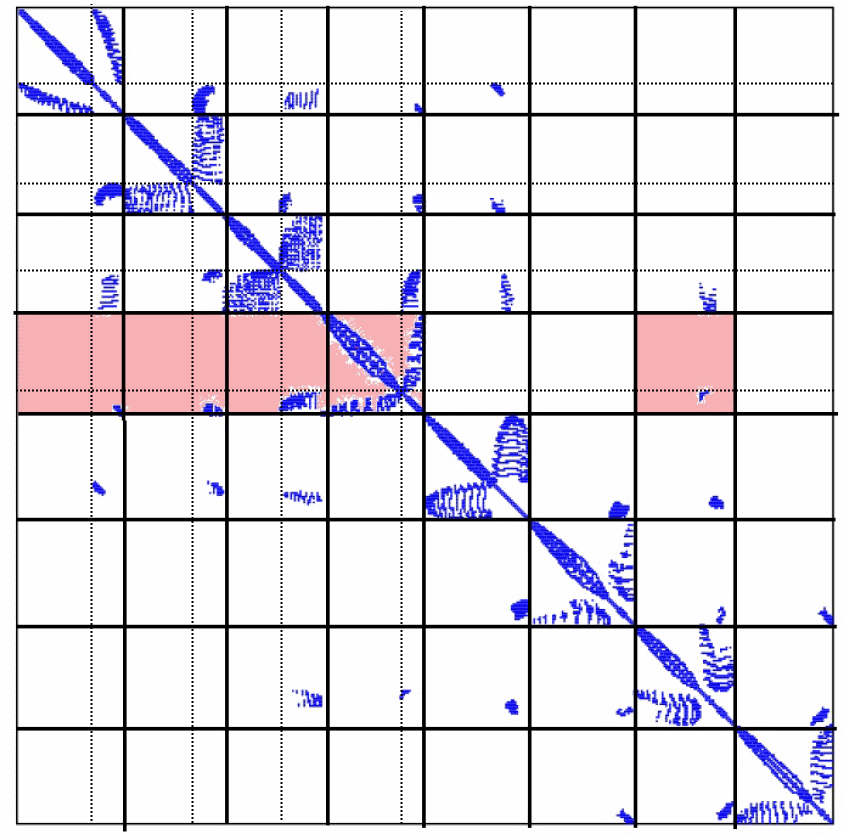}
\caption{Schematic representation of the DSMat matrix storage scheme
         implemented in Chronos using 8 MPI processes.
         The red colored blocks are assigned to process 3.\label{fig-dsmat}}
\end{figure*}

The DSMat storage \\ scheme implemented in Chronos consists in partitioning
the matrix into $n_p$ horizontal stripes of consecutive rows.
Each stripe is then divided into blocks by applying the same subdivision to the
columns, as schematically shown in Figure~\ref{fig-dsmat}, and each block is
stored as a CSR matrix.

The CSR matrices have a local numbering, i.e., rows and columns of block $IJ$
are numbered from 0 to $n_{I-1}$ and from 0 to $n_{J-1}$,
where $n_I$ and $n_J$ are the number of rows assigned to processes $I$ and $J$, respectively.
This expedient allows to use a 4-byte representation of integers, saving memory and increasing
efficiency.

Each process stores only the diagonal block and the list of "Left" (with a lower index)
and "Right" (with a higher index) blocks corresponding to the connections with
neighboring processes.
With reference to Figure~\ref{fig-dsmat}, for instance, processor 3 stores the 5 blocks
highlighted in red: 0, 1 and 2 as left neighbors, the diagonal block representing only
internal connections and 6 as right neighbor.

This blocked scheme, although a bit cumbersome to implement, allows to stress non-blocking
send/receive communications with a large superposition between communication and computation.
It has proven to be very effective in all basic operations involving a DSMat: SpMV product,
matrix-by-matrix product and matrix transposition.

\section{Numerical results}\label{sec:4}

The numerical experiments have been performed using large sparse matrices arising from
the discretization of PDEs that model challenging real-world problems.
The main goals of this section are basically to show the efficiency of the novel
implementation of the Chronos package and to demonstrate its robustness and flexibility
in dealing with severely ill-conditioned linear systems deriving form very different 
application fields, by switching the solution strategy through an appropriate parameter
tuning. As described in the previous sections, the user can exploit the effectiveness of
an advanced and tunable smoother like aFSAI, vary the coarsening ratio or switch
between different interpolation methods depending on the specific problem at hand.

Chronos is benchmarked on a set of problems that can be grouped into two classes
denoted as Fluid dynamic (F) and Mechanical (M).
The first class of benchmarks consists of a series of problems arising from the
discretization of the Laplace operator and related to fluid dynamic problems,
such as underground fluid flow (reservoir), compressible or incompressible airflow
around complex geometries (CFD) or porous flow (porous flow).
The second category includes problems related to mechanical applications such as
subsidence analysis, hydrocarbon recovery, gas storage (geomechanics),
mesoscale simulation of composite materials (mesoscale), mechanical deformation of
human tissues or organs subjected to medical interventions (biomedicine),
design and analysis of  mechanical elements, e.g., cutters, gears, air-coolers
(mechanical).

In our experiments, we consider challenging test cases, not only for the high number of
degrees of freedom (DOFs), but also because of their intrinsic ill-conditioning. Indeed, in real
applications, we usually have to deal with severe jump of the physical proprieties, complicate
geometries leading to highly distorted elements, heterogeneity and anisotropy.
The matrices considered in the experiments are listed in Table \ref{tab.matrices} with details
about the size, the number of non-zeros and the application field they arise from.
\begin{table}
\begin{center}
\centering
{\small 
\begin{tabular}{lrrrrr}
\toprule
Matrix & Class & $n$ & $nnz$ & avg. $nnz$/row & Application field \\
\midrule
{\tt finger4m}  & F &   4,718,592  &     23,591,424 &  5.00 & porous flow \\
{\tt guenda11m} & M &  11,452,398  &    512,484,300 & 44.75 & geomechanics \\
{\tt agg14m}    & M &  14,106,408  &    633,142,730 & 44.88 & mesoscale \\
{\tt M20}       & M &  20,056,050  &  1,634,926,088 & 81.52 & mechanical \\
{\tt tripod24m} & M &  24,186,993  &  1,111,751,217 & 45.96 & mechanical \\
{\tt rtanis44m} & F &  44,798,517  &    747,633,815 & 16.69 & porous flow \\
{\tt geo61m}    & M &  61,813,395  &  4,966,380,225 & 80.34 & geomechanics \\
{\tt poi65m}    & F &  65,939,264  &    460,595,552 &  6.99 & CFD \\
{\tt Pflow73m}  & F &  73,623,733  &  2,201,828,891 & 29.91 & reservoir \\
{\tt c4zz134m}  & M & 134,395,551  & 10,806,265,323 & 80.41 & biomedicine \\
\bottomrule
\end{tabular} }
\end{center}
\caption{Benchmark matrices used in the numerical experiments. For each matrix, the class,
         the size, $\mbox{n}$, the number of non-zeros, $\mbox{nnz}$, the average
         number of non-zeroes per row and the application field are provided.}
\label{tab.matrices}
\end{table}
The reader can refer to Appendix~\ref{appendix} for a detailed description of each
test case.

We subdivide the discussion of the results into two parts, the former collecting test cases from
fluid dynamics and the latter from mechanics. We also provide strong and weak scalability
analysis of the proposed implementation using large scale computational resources.
The result are presented in terms of total number of computational cores used $n_{cr}$,
the grid and operator complexities, $C_{gd}$ and $C_{op}$, respectively, the number of
iteration to converge, $n_{it}$ and the set-up, iteration and total times, $T_p$, $T_s$ and
$T_t=T_p+T_s$, respectively.\\
The right-hand side vector used for all test cases is a random vector.   
The linear systems are solved by the preconditioned conjugate gradient (PCG) method
with a zero initial solution and convergence is considered achieved when the
$l2$-norm of the iterative residual becomes smaller than $10^{-8} \cdot \|b\|$.
The Chronos performance has been evaluated on the Marconi100 supercomputer,
from the Italian consortium for supercomputing (CINECA). Marconi100,
classified within the first ten positions of the TOP500 ranking~\cite{top500} at the
time of writing, is composed by 980 nodes based on the IBM Power9 architecture, each
equipped with two 16-cores IBM POWER9 AC922 at 3.1 GHz processors.
For each test, the number of cores, $n_{cr}$, is selected to have a per core load
of about 100\--150,000 unknowns and, consequently, different numbers of nodes are allocated
for different problem dimensions.
For all the tests, each node reserved for the run is always fully exploited. Furthermore,
in order to take advantage from the hybrid implementation and use also shared memory
parallelism, we always use 8 MPI tasks on each node and 4 OpenMP threads for each task.

As a reference point to evaluate the performance of Chronos, we compare it with
the state-of-the-art solvers available from PETSc \cite{petsc-web-page}. More
specifically, we use BoomerAMG~\cite{HenYan02} and GAMG, the native PETSc aggregation-based
AMG, as preconditioners in fluid dynamics and mechanical problems, respectively.
The choice of BoomerAMG and GAMG as baseline solvers is because they are very well known
and open-source packages whose performance have been demonstrated in several papers
\cite{BreTonBec06,FalYan02,HenYan02,petsc-user-ref}.

\subsection{Fluid dynamics test cases}

The general purpose AMG implemented in Chronos is highly tunable offering
several set-up options to effectively solve this set of problems as it will be shown below.

\begin{figure}[htb]
\begin{center}
\includegraphics[width=1.0\textwidth]{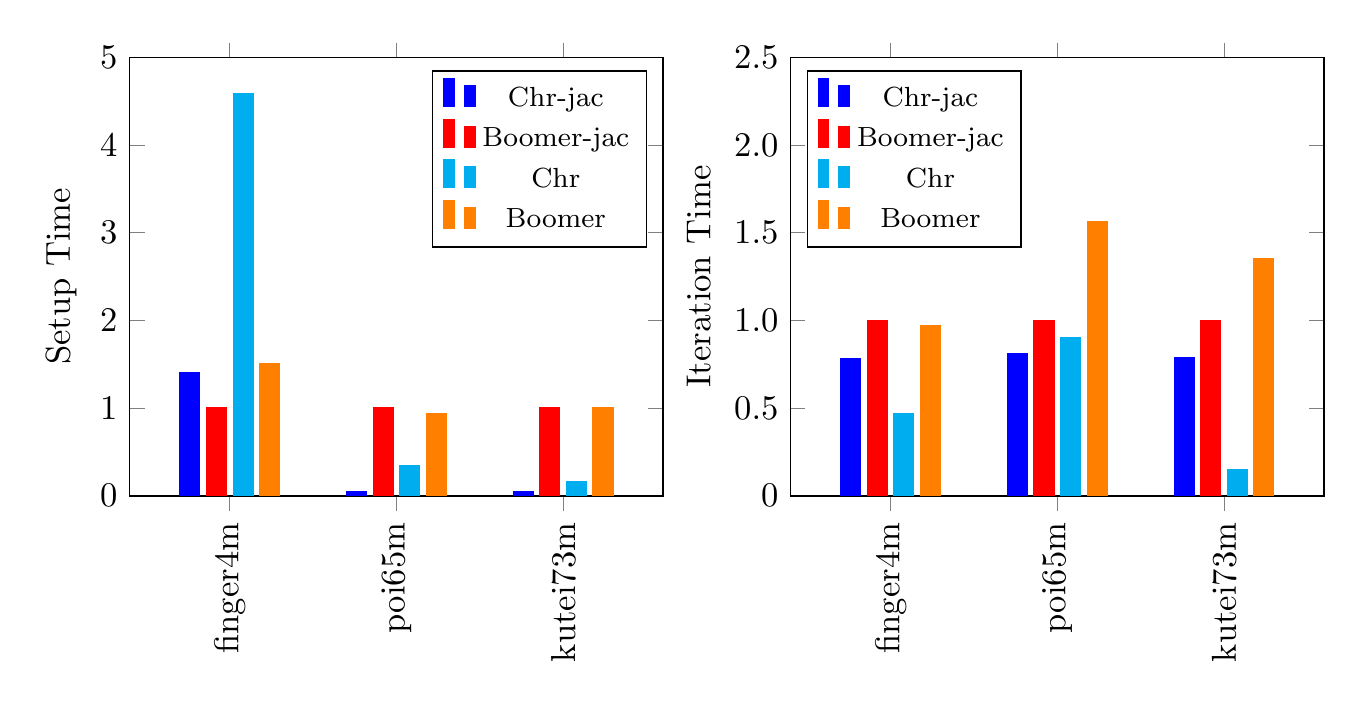}
\caption{Comparison between Chronos and BoomerAMG by using the extended+i 
prolongation and Jacobi or default smoothing.
Left: setup time. Right: solution time.}\label{fig.PETSc-aAMG}
\end{center}
\end{figure}

First, we start by comparing Chronos and BoomerAMG performance using as much as possible
the same setup. Such comparison is intended to validate our HPC implementation and to
demonstrate the efficiency of the DSMat storage scheme for SpMV product.
To this purpose, we consider the three test cases {\tt finger4m}, {\tt poi65m} and {\tt Pflow73m}.
The comparison takes place with the same preconditioner configuration, i.e.,
Jacobi smoothing, classical SoC with $\theta = 0.25$, PMIS coarsening and extended+i
prolongation.
%
The first two rows of each test case reported in Table \ref{tab.results_F} provide the
results obtained with this {\em standard set-up}.
We denote by Chr-jac and Boomer-jac the Chronos and BoomerAMG preconditioners, both
paired with Jacobi smoothing.
\begin{table}
\begin{center}
\centering
{\small 
\begin{tabular}{lrrrrrrrr}
\toprule
Matrix & $n_{cr}$ & Solv.\@ type & $C_{gd}$ & $C_{op}$ & $n_{it}$ & $T_p$ [s] &
$T_s$ [s] & $T_t$ [s] \\
\midrule
{\tt finger4m}  &  32 & Chr-jac    & 1.453 & 2.558 &  16 &  1.13 &  0.55 &  1.68 \\
                &  32 & Boomer-jac & 1.454 & 2.574 &  16 &  0.81 &  0.70 &  1.51 \\
                &  32 & Chr        & 1.453 & 2.558 &   7 &  3.71 &  0.33 &  4.04 \\
                &  32 & Boomer     & 1.454 & 2.574 &  12 &  0.79 &  0.94 &  1.73 \\
\midrule
{\tt poi65m}    & 384 & Chr-jac    & 1.327 & 4.036 &  16 &  3.81 &  1.65 &  4.46 \\
                & 384 & Boomer-jac & 1.361 & 4.450 &  13 &  84.6 &  2.03 &  86.7 \\
                & 384 & Chr        & 1.346 & 4.496 &   6 &  27.5 &  1.84 &  29.34 \\
                & 384 & Boomer     & 1.361 & 4.450 &  14 &  80.2 &  3.18 &  83.4 \\
\midrule
{\tt Pflow73m}  & 480 & Chr-jac    & 1.125 & 1.614 & 3308 &   14.1 &  611.9  &  626.0 \\
                & 480 & Boomer-jac & 1.123 & 1.593 & 3576 &  336.5 &  771.7  & 1108.2 \\
                & 480 & Chr        & 1.123 & 2.346 &  410 &   57.7 &  120.9  &  178.6 \\
                & 480 & Boomer     & 1.123 & 1.593 & 2777 &  340.5 & 1042.3  & 1382.8 \\
\bottomrule
\end{tabular} }
\end{center}
\caption{Solution of three fluid dynamic test cases among those reported in
Table~\ref{tab.matrices}. For each run, the following information is provided:
the number of cores $n_{cr}$, the grid $C_{gd}$ and operator $C_{op}$ complexities,
the number of PCG iteration $n_{it}$, the set-up time $T_p$, the iteration time $T_s$
and the total time $T_t$.}\label{tab.results_F}
\end{table}
First, we observe that the grid and operator complexities obtained with the two software
are basically the same and also the iteration count turns out to be quite similar, showing
that the two implementations are consistent. Only a slight difference occurs for
{\tt Pflow73m} but we believe it is compatible with very small differences in the code
implementations.\\ 
Figure \ref{fig.PETSc-aAMG} provides the time spent for the preconditioner set-up (left)
and for the conjugate gradient iterations (right) for each solving strategy.
Each time reported in the figure is normalized with respect to the value obtained with
Boomer-jac, which is our baseline.
We can observe that Chronos is faster than BoomerAMG in the set-up for {\tt poi65m} and
{\tt Pflow73m}, while BoomerAMG is better in {\tt finger4m}. 
Differently, as far as the solving time is concerned, Chronos slightly outperforms
BoomerAMG in all the tests.
In all the cases tested, the SpMV and the Chronos implementation turns out to be
very efficient and the total solution time obtained is comparable and sometimes
even much better than those obtained with the BoomerAMG thanks to a faster set-up.

In Table~\ref{tab.results_F} again, on the third and fourth rows, the labels Chr and
Boomer identify the results obtained with Chronos and BoomerAMG when the default
smoothers are selected, i.e., aFSAI and hybrid Gauss-Seidel, respectively. 
The use of a more evolute smoother with respect to either Jacobi or hybrid Gauss-Seidel
gives a significant advantage in terms of iteration count and solving time at the price
of a more expensive set-up, as shown in Table~\ref{tab.results_F} and
Figure~\ref{fig.PETSc-aAMG}.
The use of the aFSAI always allows for achieving a faster convergence. Furthermore,
the more ill-conditioned the problem is, the better aFSAI compares with other smoothers.
In {\tt Pflow73m}, which is the hardest problem in fluid dynamics, Chronos with aFSAI
smoothing is 6 times faster than BoomerAMG. The set-up time is larger, but
the speed-up obtained in the iteration stage may justify this effort, especially
in transient simulations where the user may have to solve repeatedly the same linear
system and can take advantage of preconditioner recycling.

In fluid dynamics, the prolongations of choice in combination with classical AMG are typically
the classical or extended+i interpolations. This last is usually more effective,
although more expensive, for changeling problems due to its ability to accurately
interpolate also fine nodes having strong fine neighbors that do not share the same
strong coarse node, possibly produced by high coarsening ratios.
In Table \ref{tab.results_prol_F}, we compare these two well-known prolongations to
the hybrid one that has been widely discussed in section~\ref{subsec:Int}.
\begin{table}
\begin{center}
\centering
{\small 
\begin{tabular}{lrrrrrrrr}
\toprule
Matrix & $n_{cr}$ & Prol.\@ type & $C_{gd}$ & $C_{op}$ & $n_{it}$ & $T_p$ [s] & $T_s$ [s] & $T_t$ [s] \\
\midrule
                &  32 & Chr-clas    & 1.467 & 1.871 &  31 &  3.57 &  1.32 &  4.89 \\
{\tt finger4m}  &  32 & Chr-hybc    & 1.465 & 2.051 &  14 &  3.62 &  0.66 &  4.28 \\
                &  32 & Chr-exti    & 1.453 & 2.558 &   7 &  3.71 &  0.33 &  4.04 \\
\midrule	    						
                & 384 & Chr-clas    & 1.612 & 1.943 &  46 &  23.3 &  6.47 &  29.8 \\
{\tt rtanis44m} & 384 & Chr-hybc    & 1.585 & 2.030 &  36 &  26.6 &  6.09 &  32.8 \\
                & 384 & Chr-exti    & 1.572 & 2.580 &  16 &  34.0 &  2.90 &  36.9 \\
\midrule
                & 384 & Chr-clas    & 1.381 & 2.339 &  21 &  17.9 &  4.69 &  22.59 \\
{\tt poi65m}    & 384 & Chr-hybc    & 1.361 & 2.888 &  13 &  19.0 &  2.57 &  21.57 \\
                & 384 & Chr-exti    & 1.346 & 4.496 &   6 &  27.5 &  1.84 &  29.34 \\
\midrule
                & 480 & Chr-clas    & 1.236 & 1.391 & 414 &  39.4 & 60.2  & 99.6 \\
{\tt Pflow73m}  & 480 & Chr-hybc    & 1.234 & 1.448 & 416 &  40.4 & 67.1  & 107.5 \\
                & 480 & Chr-exti    & 1.234 & 2.346 & 410 &  57.7 & 120.9 & 178.6 \\
\bottomrule
\end{tabular} }
\end{center}
\caption{Comparison between different interpolation formulas in the
solution of the fluid dynamic test problems from Table \ref{tab.matrices}.
For each run, the following information is provided:
number of cores $n_{cr}$, prolongation type, grid $C_{gd}$ and operator $C_{op}$
complexities, number of iteration $n_{it}$, set-up time $T_p$, iteration time $T_s$ and
total time $T_t$.}
\label{tab.results_prol_F}
\end{table}
Let us consider first the results obtained for {\tt finger4m} and {\tt poi65m} for which the
solver behavior is quite similar. We can observe that the extended+i interpolation
is the more accurate one, with the higher value for operator complexity.
As expected, this leads to a lower number of iterations, but a higher computational
cost per iteration.
On the contrary, the classical interpolation formula is the cheapest to compute, 
with a very low operator complexity. However, taking into account only distance-one coarse
nodes, the prolongation operator is not able to efficiently reproduce the smooth error, 
causing an increase of the iteration count, up to twice the iteration count obtained with
extended+i. For these two tests, the best configuration lies in the
middle of these two, i.e., the hybrid interpolation formula, which keeps low the operator
complexity taking into account just the distance-two coarse nodes actually useful to the
interpolation process.  In this way, we are able to obtain a more accurate interpolation
formula with a computational cost comparable to the classical one.

The behavior is quite different for the other two test cases.\\
In {\tt rtanis44m}, we have a strong heterogeneity and anisotropy of the permeability tensor,
factors that dramatically increase the problem ill-conditioning.
Hence, the most accurate interpolation method, i.e., extended+i, is needed to efficiently
solve this problem. The iteration count is one third with respect to classical interpolation
and the solution time is approximately one half. 
Unlike before, the increased accuracy of the hybrid interpolation over classical is not enough
to give a sufficient benefit in terms of solving time.
It is worth noting that the increased set-up cost for extended+i is in this case largely
compensated in the iteration stage. This gain is even more pronounced in cases where
preconditioner recycling is possible such as in some transient or non-linear simulations.\\
The last test case considered in this section is {\tt Pflow73m}, a very challenging and
severely ill-conditioned problem from underground flow.
Even if this is a diffusion problem, the great jumps in permeability and the distorted mesh
lead to a matrix whose near-kernel is not well represented by the unitary vector.  
For this reason, the number of iterations required to achieve the convergence increases a lot
with respect to the other tests and not even the most accurate interpolations
such as extended+i or hybrid give any benefit over classical interpolation.
Hence, the cheapest classical formula proves also the most efficient strategy for this
test case. 

\subsection{Mechanical test cases}

In this section, the potential of Chronos and its effectiveness in mechanical problems
are highlighted. \\
As seen above, Chronos allows for setting-up a very flexible AMG preconditioner,
adaptable to problem types the user has to solve, with different choices available
for interpolation operators and smoothing methods.
In addition, it allows the possibility to directly smooth the prolongation
and/or filter it. 
As in the previous paragraph, we first define a baseline with
state-of-the-art methods such as BoomerAMG (Boomer), with Hybrid Gauss-Seidel smoothing,
the unknown-based Boomer with separate treatment of unknowns relative to different
directions (unk-based-Boomer) and the GAMG, an aggregation-based method.\\
We first refer to the test case {\tt tripod24m}, whose results are provided in Table
\ref{tab.results_methods}.
With the standard Boomer, the solution is reached with a high number of iterations,
more than 900 and the iteration time responsible of most of the total solution time.
A significant improvement is obtained using the unknown-based version~\cite{BakKolYan09},
where iterations are reduced by one third, and set-up and iteration times drop by 50\%.
The aggregation based AMG seems to be the most effective one for mechanical problems as,
with GAMG, iterations are further reduced, and both $T_p$ and $T_s$ times decrease
significantly. In this problem, Chronos with BAMG prolongation and aFSAI smoother
(BAMG-aFSAI) is more effective than GAMG with a speed-up of two on the total time.
The set-up time is larger, but the number of PCG iterations is lower and the cost
per iteration is one-third of that of GAMG.
It is also possible to smooth the prolongation operator with Jacobi. We denote this
method as SBAMG-aFSAI. As could be expected, the operator complexity and the set-up
time both increase but, on the other hand, the number of iterations to converge and
the solution time are smaller.
Operator complexity and set-up time increases can be limited by means of filtering
(FBAMG-aFSAI) without compromising effectiveness. FBAMG-aFSAI requires the same number
of iterations to converge but at a lower cost per iteration.
These two last strategies are particularly effective in a FEM simulation where
the preconditioner can be reused several times in different time-steps so that the
set-up cost becomes secondary.

\begin{table}
\begin{center}
\centering
{\small 
\begin{tabular}{lrrrrrrrr}
\toprule
Matrix & $n_{cr}$ & Prec.\@ type & $C_{gd}$ & $C_{op}$ & $n_{it}$ & $T_p$ [s] & $T_s$ [s] & $T_t$ [s] \\
\midrule

                &      & Boomer           & 1.244 & 3.207 &   931 &  64.1 & 611.9 & 676.1\\
                &      & unk-based-Boomer & 1.328 & 3.669 &   335 &  43.8 & 262.4 & 203.5\\
{\tt tripod24m} &  160 & GAMG             & 1.543 & -     &   294 &  12.1 &  80.5 &  92.6\\
                &      & BAMG-aFSAI       & 1.041 & 1.116 &   222 &  21.8 &  23.0 &  44.8\\
                &      & SBAMG-aFSAI      & 1.041 & 1.322 &   118 &  36.7 &  16.1 &  52.9\\
                &      & FBAMG-aFSAI      & 1.041 & 1.212 &   120 &  33.5 &  13.5 &  47.0\\

\bottomrule
\end{tabular} }
\end{center}
\caption{Solution of the {\tt tripod24m} test case from Table \ref{tab.matrices} with
different approaches.
For each run, the following information is provided:
the number of cores $n_{cr}$, the preconditioner type, the grid complexity $C_{gd}$,
the operator complexity $C_{op}$, the number of iteration $n_{it}$,
the set-up time $T_p$, the iteration time $T_s$, and the total time $T_t$.}\label{tab.results_methods}
\end{table}

\begin{figure}[htb]
\begin{center}
\includegraphics[width=1.0\textwidth]{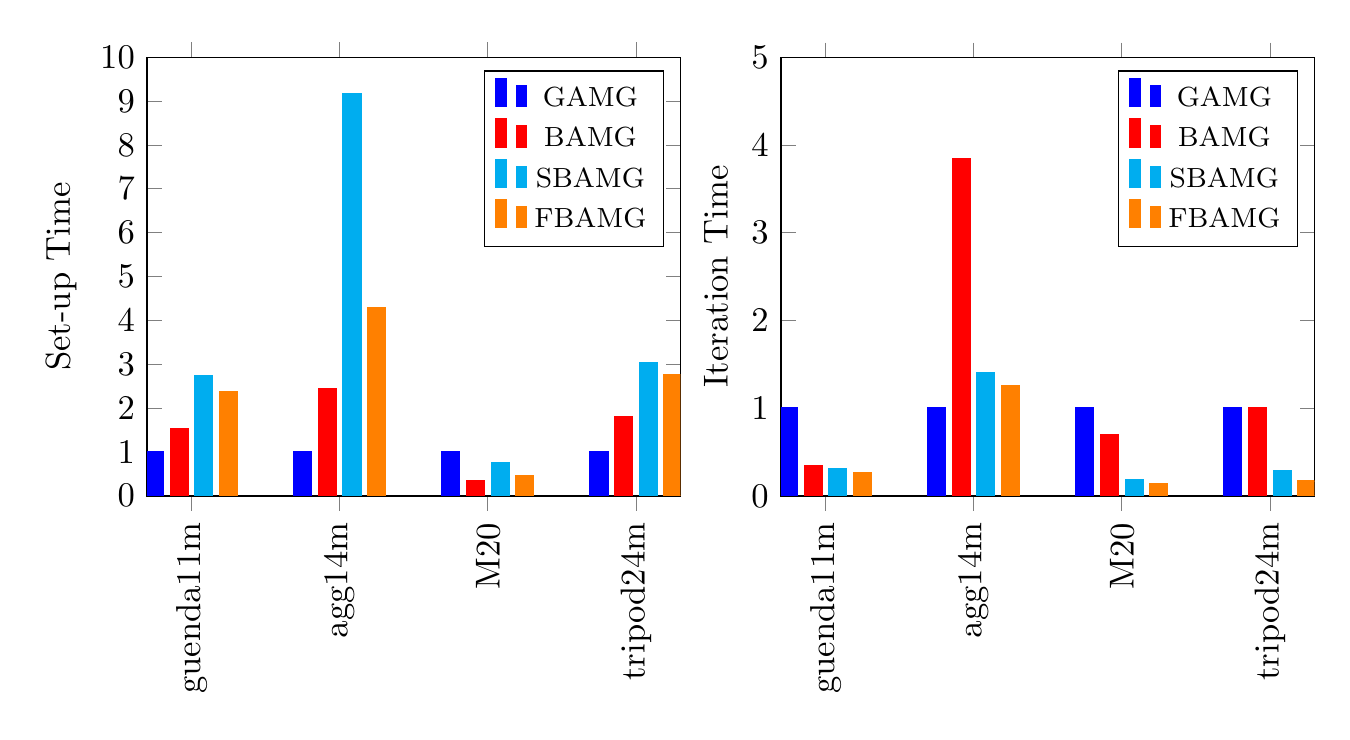}
\caption{Comparison between GAMG and the BAMG strategies on the mechanical test cases.
Left: normalized $T_p$ to the GAMG solution. Right: normalized $T_s$ to the GAMG solution.}
\label{fig.GAMG-BAMG}
\end{center}
\end{figure}

Chronos proved robust and efficient in addressing all the mechanical test cases.
A comparison of the number of iterations and times obtained with GAMG and the three BAMG
strategies outlined above is shown in Table \ref{tab.GAMG-BAMG}.
To highlight the speed-up, Figure \ref{fig.GAMG-BAMG} shows set-up and the iteration time
normalized to the GAMG times.
Unfortunately, the comparison for the two largest cases, {\tt geo61m } and {\tt c4zz13m},
is not reported because these matrices have not been dumped on file due to their
large size, and the tests have been run by linking Chronos to the FEM program
ATLAS~\cite{ATLAS-webpage}.
For the three benchmarks, {\tt guenda11m}, {\tt tripod24m} and {\tt M20}, the number
of PCG iterations required by GAMG and BAMG is comparable, but the overall solution time
is significantly lower for BAMG with a speed-up of Chronos over GAMG up to 4 in
these tests. The only exception is the matrix {\tt agg14m} where GAMG is able to produce
a very effective preconditioner at the lowest set-up cost.

\begin{figure}[htb]
\begin{center}
\includegraphics[width=0.7\textwidth]{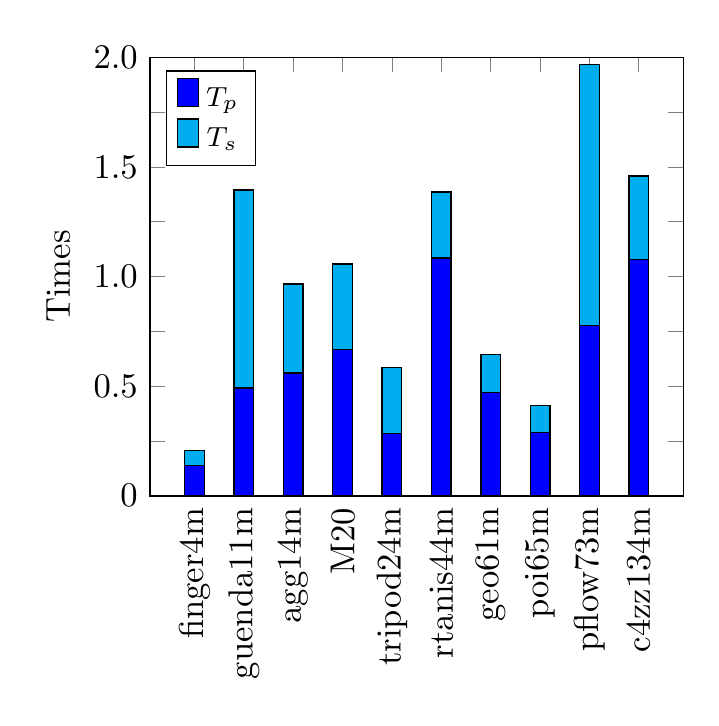}
\caption{Set-up and iteration times, for all benchmark problems, normalized over
         the resources allocated per non-zero.}
\label{fig.NORM}
\end{center}
\end{figure}

Finally, with the aid of Figure~\ref{fig.NORM}, we would like to point out how the total
solution time depends only mildly on the problem nature but on its size only.
Figure~\ref{fig.NORM} shows for each problem the total solution time divided by the number
of non-zeroes per allocated core, and this resulting time is further normalized with the
average among all the experiments. In other words, the figure should show the solution
time for each problem as if {\em exactly} the same resources were allocated for each
non-zero. For a preconditioner that is totally independent by the problem nature,
it would be expected the same solution time for every problem. It can be observed that,
through careful parameter tuning, Chronos is able to produce total solution times
very close to the average normalized solution time, thus showing only a mild dependence
on the application at hand.

\begin{table}
\begin{center}
\centering
{\small 
\begin{tabular}{lrrrrrrrr}
\toprule
Matrix & $n_{cr}$ & Prol.\@ type & $C_{gd}$ & $C_{op}$ & $n_{it}$ & $T_p$ [s] & $T_s$ [s] & $T_t$ [s] \\
\midrule
                 &   64 &  GAMG	 &  1.580	& -	    & 978	&  18.3 &  306.2 &  324.5\\
{\tt guenda11m}  &   64 &  BAMG	 &  1.041	& 1.118	& 937	&  27.8 &  105.0 &  133.0\\
                 &   64 &  SBAMG &  1.041	& 1.354	& 638	&  50.3 &   96.3 &  147.0\\
                 &   64 &  FBAMG &  1.041	& 1.240	& 638	&  43.5 &   79.8 &  123.0\\
\midrule
                 &  128 &  GAMG	 &  1.644	& -     &  26   &  12.5 &    5.8 &   18.2\\
{\tt agg14m}     &  128 &  BAMG	 &  1.085	& 1.287	& 135	&  30.6 &   22.2 &   52.8\\
                 &  128 &  SBAMG &  1.085	& 2.264	&  31	& 114.4 &    8.1 &  122.6\\
                 &  128 &  FBAMG &  1.085	& 1.670	&  34	&  53.6 &    7.3 &   60.9\\
\midrule
                 &  128 &  GAMG	 &  1.162 & -	    & 245	& 211.0 &  391.4 &  602.4\\
{\tt M20}        &  128 &  BAMG	 &  1.054	& 1.184	& 775	&  71.2 &  275.0 &  347.0\\
                 &  128 &  SBAMG &  1.054	& 1.677	& 151	& 158.0 &   71.2 &  229.2\\
                 &  128 &  FBAMG &  1.054	& 1.292	& 158	&  93.9 &   55.1 &  149.1\\
\midrule
                 &  160 &  GAMG	 &  1.543	& -	    & 294	&  12.1 &   80.5 &   92.6\\
{\tt tripod24m}  &  160 &  BAMG	 &  1.041	& 1.116	& 222	&  21.8 &   23.0 &   44.8\\
                 &  160 &  SBAMG &  1.041	& 1.322	& 118	&  36.7 &   16.1 &   52.9\\
                 &  160 &  FBAMG &  1.041	& 1.212	& 120	&  33.5 &   13.5 &   47.0\\
\bottomrule
\end{tabular} }
\end{center}
\caption{Comparison between different interpolation formula in the
solution of the mechanical test problems from Table \ref{tab.matrices}.
For each run, the following information is provided:
number of cores $n_{cr}$, prolongation type, grid $C_{gd}$ and operator $C_{op}$
complexities, number of iteration $n_{it}$, set-up time $T_p$, iteration time $T_s$ and
total time $T_t$.}
\label{tab.GAMG-BAMG}
\end{table}

\subsection{Strong and weak scalability}

\begin{figure}[htb]
\begin{center}
\includegraphics[width=0.5\textwidth]{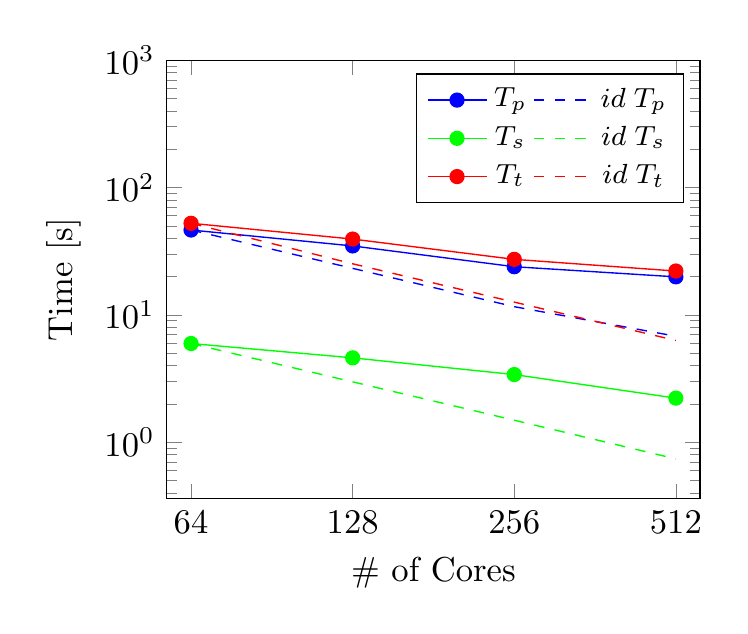}
\includegraphics[width=0.49\textwidth]{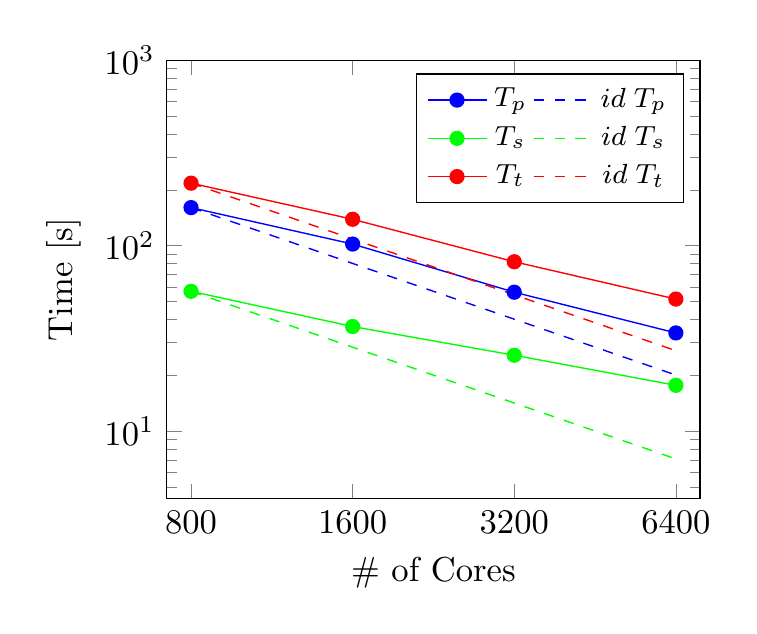}
\caption{Strong scalability test for {\tt poi65m} matrix and BAMG prolongation (left)
         and {\tt c4zz134m} matrix and extended+i prolongation (right).
         Preconditioner set-up time $T_p$, iteration time $T_s$ and total time $T_t$
         are provided.}
\label{fig.strong}
\end{center}
\end{figure}

\begin{figure}[htb]
\begin{center}
\includegraphics[width=1.0\textwidth]{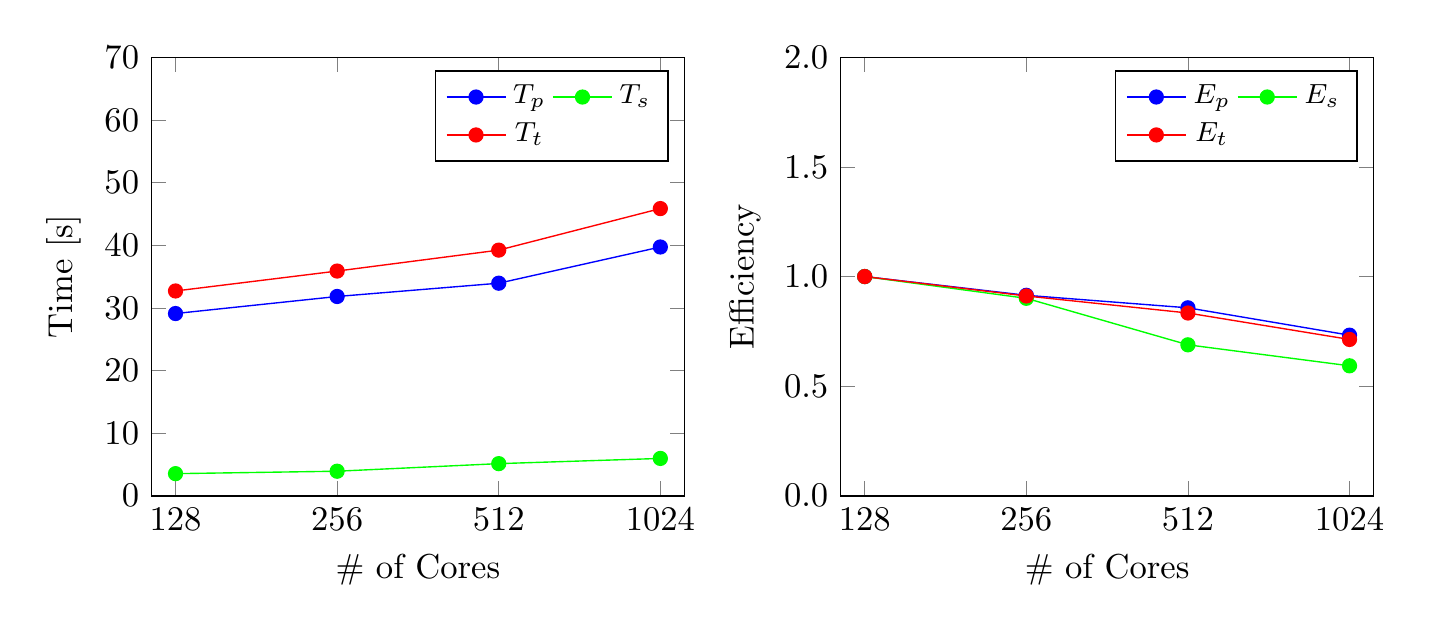}
\caption{Week scalability test on a standard 7-point finite difference discretization of
         the Poisson problem. Set-up, iteration and total times vs. $n_{cr}$ are reported on
         the left, corresponding efficiencies vs. $n_{cr}$ on the right.}
\label{fig.week_Poi}
\end{center}
\end{figure}

In this last subsection, we evaluate the strong and weak scalability of the
AMG preconditioners implemented in Chronos.
All the three times, i.~e. set-up $T_p$, iteration $T_s$ and total $T_t$ times, 
are analyzed to assess scalability.
The strong scalability test is shown in Figure \ref{fig.strong}, on the left
for the {\tt c4zz134m} test matrix with BAMG prolongation and on the right for {\tt poi65m}
with extended+i prolongation.
The number of cores varies from the minimum necessary to store matrix and preconditioner
up to 8 times the initial number. In both tests, the times decrease as the computing
resources increase, with a trend close to the ideal one.\\ 
Finally, the weak scaling is investigated with a standard 7-point finite difference
discretization of the Poisson problem.
Figure \ref{fig.week_Poi} shows, on the left, both the total time spent in the set-up and solve
phase and, on the right, the parallel efficiency.
Efficiency of weakly scaling up to $N$ nodes is defined as $E = T_N/(N T_1)$, with $T_1$
the time required on a single node and $T_N$ the time on $N$ nodes. 
In this test, we always assign 218,750 unknowns per core.

The result shows that efficiency is very good and keeps almost constant in the first two
doubles of the cores, whereas a bit greater efficiency drop occurs in the last one. This
performance dropdown can be ascribed to two different factors. First of all, while
Marconi100 cores can be fully reserved for the test runs, the overall network is always
shared with other users, and, consequently, the larger is the resource allocation,
the larger the disturbance from other running processes. Secondly, a performance dropdown
is almost unavoidable in AMG methods, as the grid hierarchy ends always up with small grids.
The larger the number of resources allocated, the less efficient will be the software
in dealing lower levels. Currently, to ease the implementation, Chronos uses all the
allocated cores on each grid except the last one, where an {\tt allgather} operation
is called from a single core to solve the coarsest problem. In a future implementation,
we plan to progressively reduce the amount of resources with levels, thus reducing the
network traffic and increasing efficiency.

\section{Conclusions}
\label{sec:5}

In this work, the Chronos library for the solution of large and sparse linear
algebra problems on high performance platforms has been presented with a deep
analysis of its numerical and computational performance.
Chronos, which will be freely accessible to research
institutions~\cite{CHRONOS-webpage}, provides iterative solution methods for
linear systems and eigenproblems along with advanced parallel preconditioners.

Although the library comprises classical and novel methods already known in the literature,
all of its algorithms have been attentively revisited, tuned and optimized on the basis
of a large experimentation on real-world and industrial benchmarks arising from very different
application fields. Moreover, every numerical kernel has been designed with special
attention to its parallel performance and future extensibility to new
numerical approaches and hardware.

The wide set of numerical experiments, provided in the work, clearly shows the
ability of Chronos to give excellent performance in very different applications
with solution times no worse or even better than those offered by other widely
used HPC linear solvers as BoomerAMG and GAMG. Furthermore, this library offers
great flexibility in the choice of the preconditioning strategy with the result
that, once a proper set-up is found, total solution time depends solely or almost
solely on the problem size and the number of computational resources allocated.

Our future work will be focused on porting Chronos on more energy-efficient and
promising hardware such as GPU accelerators or FPGA, as well as using the innermost
kernels of the library in developing advanced block preconditioners for multi-physics
applications.

We also plan to build a stronger theoretical basis for the adaptive construction
of the test space, unavoidable in problems lacking an initial guess for the operator
near kernel, and for the operator and prolongation filtering which can greatly
improve performance in though problems where denser operators may be needed.

\vspace{0.5cm}
\noindent {\bf Acknowledgments.} The authors gratefully thank Proff. S. Koric, G. Mazzucco
and E.L. Carniel, who provided the matrices M20, agg14m and c4zz134m used in the experiments.

\appendix
\section{Description of the real-world applications}
\label{appendix}

This appendix provides a detailed description of the test cases presented in this work, 
as listed in Table \ref{tab.matrices}. 

\subsection{Test case {\tt finger4m}}  
\label{sec:finger4m}

The matrix {\tt finger4m} derives from a two-dimensional Darcy flow of a binary mixture.
The physical model describes flow in a porous medium or Hele-Shaw cell, a thin gap
between two parallel plates. The behavior of the system, and hence the matrix, is governed
by two nondimensional groups: the P\'eclet number $Pe = 10^4$ and the viscosity ratio
$M=\exp(3.5)$ \cite{jha2011fluid}.

\subsection{Test cases {\tt guenda11m} and {\tt geo61m}} 
\label{sec:guenda11m}

The matrices {\tt guenda11m} and {\tt geo61m} derive from two 3D geomechanical models of
a reservoir.

\begin{figure}[!htbp]
\centering
\includegraphics[width=0.8\linewidth]{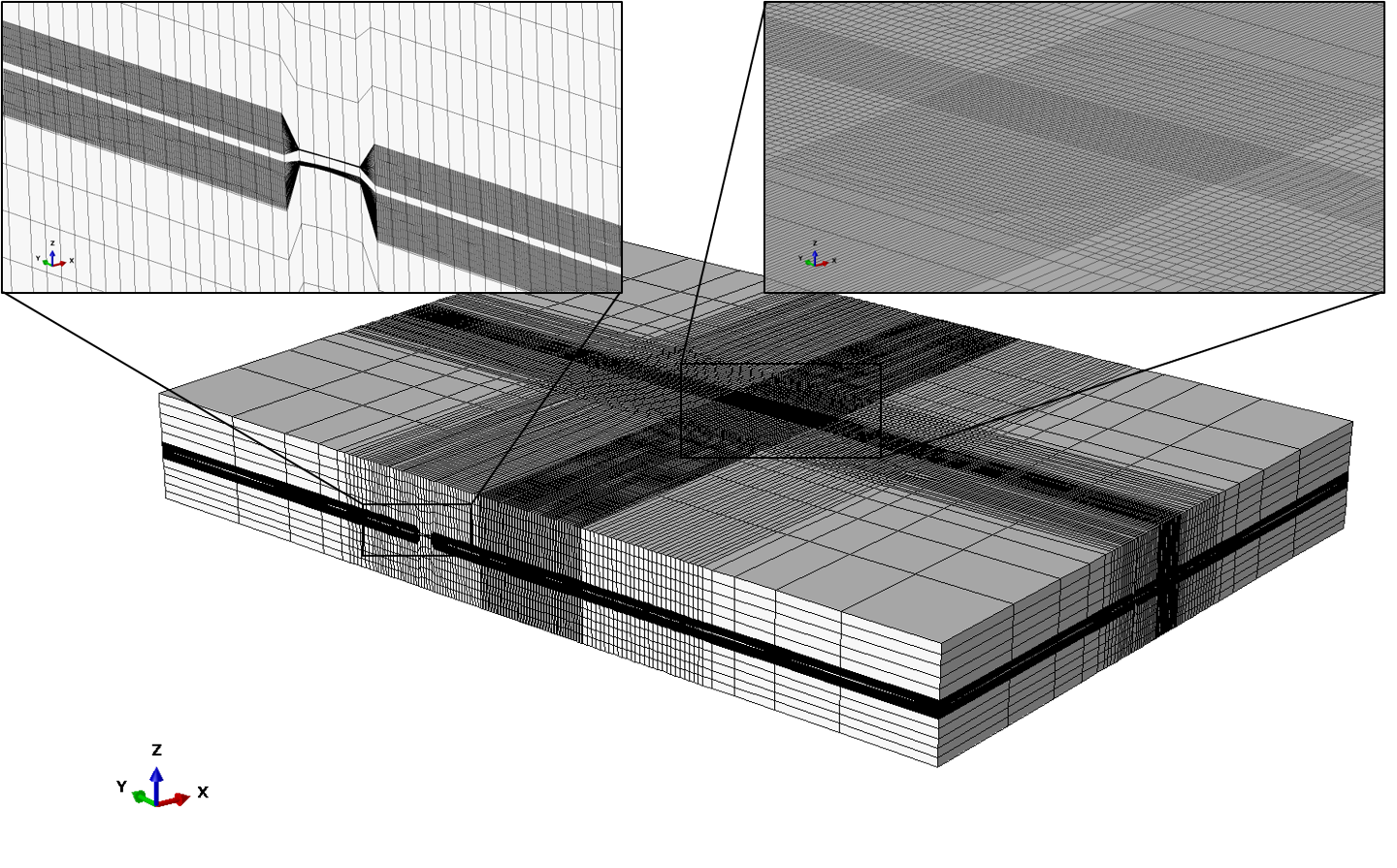}
\caption{Geometry of the test case {\tt geo61m}. The box show some are where the mesh id
         very fine and eventually with distorted elements.}
\label{fig:guenda11_mesh}
\end{figure}

In particular, the matrix {\tt guenda11m} derives from a domain that spans an area of
$40 \times 40 \, km^{2}$ and extends down to $5 \, km$ depth. To reproduce with high fidelity
the real geometry of the gas reservoir, a severely distorted mesh with  $22,665,896$ linear
tetrahedra and $3,817,466$ vertices is used. While fixed boundaries are prescribed on the
bottom and lateral sides, the surface is traction-free.\\
The matrix {\tt geo61m} represents a geological formation with 479 layers. The geometry of
the modeled domain is characterized by an area of $55 \times 40 \, km^{2} $ with the reservoir
in an almost barycentric position and the base at a depth of $6.5 \, km^{} $. The grid is
based on a mesh of $20,354,736$ brick elements. The Figure \ref{fig:guenda11_mesh} shows a
representation of the problem's geometry and mesh. As can be observed, some elements are
highly distorted to reproduce the geological layers.

\subsection{Test case {\tt agg14m}} 
\label{sec:agg14m}

\begin{figure}[!htbp]
\centering
\includegraphics[width=0.8\linewidth]{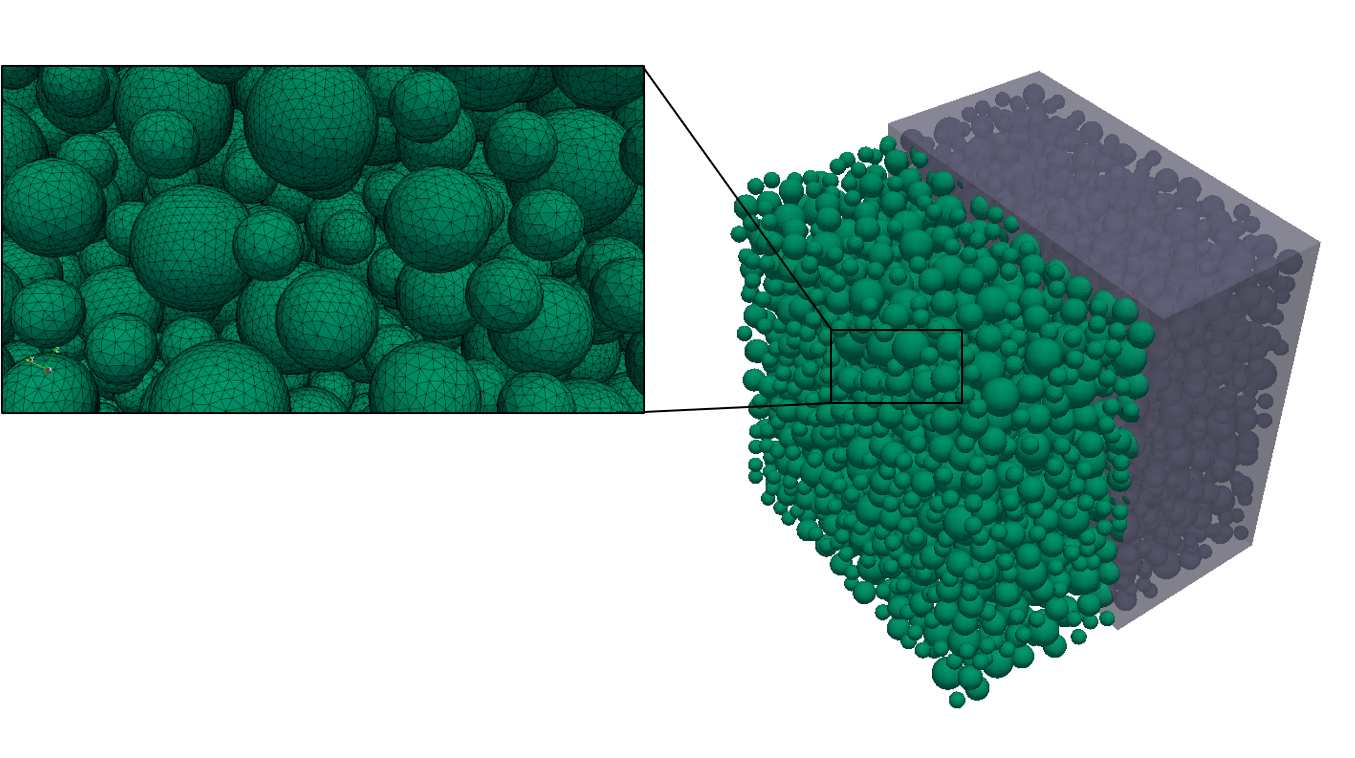}
\caption{Geometry of the test case {\tt agg14m}. The box represents the mesh of the
         spherical inclusions.}
\label{fig:agg14m_mesh}
\end{figure}

The mesh derives from a 3D mesoscale simulation of an heterogeneous cube of lightened concrete.
The domain has dimensions $50 \times 50 \times 50 \, mm^{3}$ and contains $2644$ spherical
inclusions of polystyrene. The cement matrix is characterized by
$(E_1, \nu_1) = (25,000 MPa, 0.30)$, while the polystyrene inclusions are characterized by
$(E_2, \nu_2) = (5 MPa, 0.30)$ \cite{mazzucco2018numerical,mazzucco2020tomography}. Hence,
the contrast between the Young modules of these two linear elastic materials is extremely high. 
The discretization is done via tetrahedral finite elements. The Figure \ref{fig:agg14m_mesh}
shows a representation of the problem's geometry and mesh.

\subsection{Test case {\tt M20}} 
\label{sec:M20}

The mesh derives from the 3D mechanical equilibrium of a symmetric machine cutter that
is loosely constrained. The unstructured mesh is composed by $4,577,974$ second order
tetrahedra and $6,713,144$ vertices resulting in $20,056,050$ DOFs. Material is linear
elastic with $(E,\nu) = (10^{8} MPa, 0.33)$. 
This problem was initially presented by \cite{KorLuGul14} and later used in the work
\cite{KorGup16}.

\subsection{Test case {\tt tripod24m}} 
\label{sec:tripod3239k}

\begin{figure}[!htbp]
\centering
\includegraphics[width=0.8\linewidth]{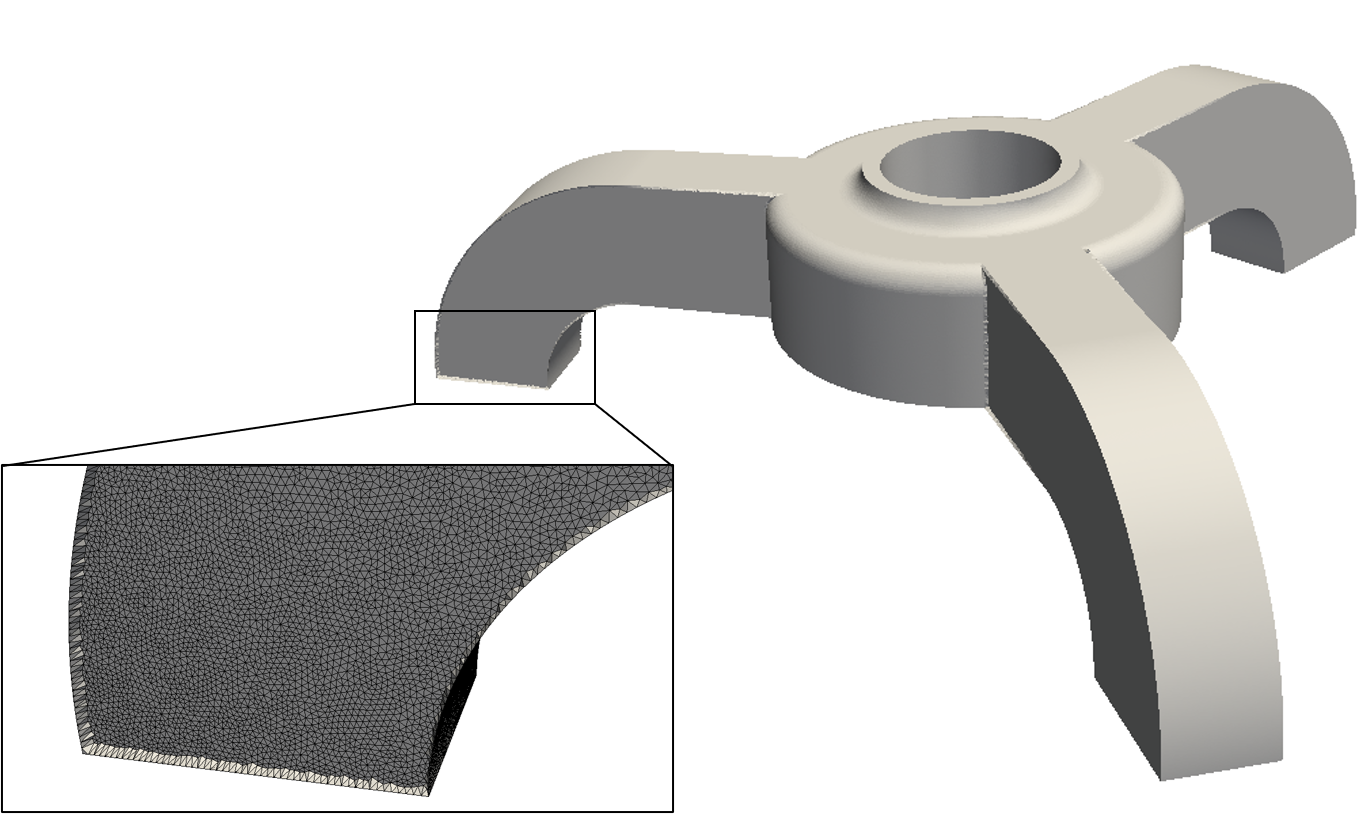}
\caption{Geometry of the test case {\tt tripod24m}. The box represents the mesh with
         (almost regular) linear tetrahedra.}
\label{fig:tripod_mesh}
\end{figure}

The mesh derives from the 3D mechanical equilibrium of a tripod with
clamped bases. Material is linear elastic with $(E,\nu) = (10^{6} MPa, 0.45)$. The mesh is
formed by linear tetrahedra and discretization is given by the finite
element method. Figure \ref{fig:tripod_mesh} shows the
geometry and the mesh of the problem.

\subsection{Test case {\tt rtanis44m}}  
\label{sec:rtanis44m}

The mesh derives from a 3D diffusion problem in a porous media.
The diffusion problem is governed by an anisotropic permeability tensor of the form
$\hat{K} =Q^{T}KQ$, where $Q$ is a rotation matrix and $K$ is a diagonal matrix defined as

\[ Q = \left( \begin{array}{ccc}
\cos(\theta) & -\sin(\theta) & 0 \\
\sin(\theta) & \cos(\theta) &  0\\ 
0 & 0 & 1
\end{array} \right), 
K = \left( \begin{array}{ccc}
K_x &0 &0 \\
0 & K_y & 0 \\
0 & 0 & K_z 
\end{array} \right),
\]
with the rotation angle $\theta = 30^{\circ}$ and the permeability matrix given by $K_x=10.0$,  $K_y =1.0^{-3}$, $K_z =1.0^{-6}$.

\subsection{Test case {\tt poi65m}} 
\label{sec:poi65m}
The mesh derives from the solution of the Poisson's equation $\nabla^2 \phi = f$  over a
3D cube. The domain is discretized with a $100 \times 200 \times 402$ finite difference
grid.

\subsection{Test case {\tt Pflow73m}} 
\label{sec:Pflow73m}

The mesh derives from a basin model, with the discretization of a
$178.8 \times 262.0 \, km^{2}$ geological area - at the end of basin evolution -  with
a mesh of 20-node hexahedral elements. The {\tt Pflow73m} matrix derives from the
discretization of the mass conservation and Darcy’s law. Due to strong permeability
contrasts between neighboring elements and geometrical distortion of the computational
grid, the matrix is severely ill-conditioned and challenging to solve.

\subsection{Test case {\tt c4zz134m}} 
\label{sec:c4zz134m}
\begin{figure}[!htbp]
\centering
\includegraphics[width=0.8\linewidth]{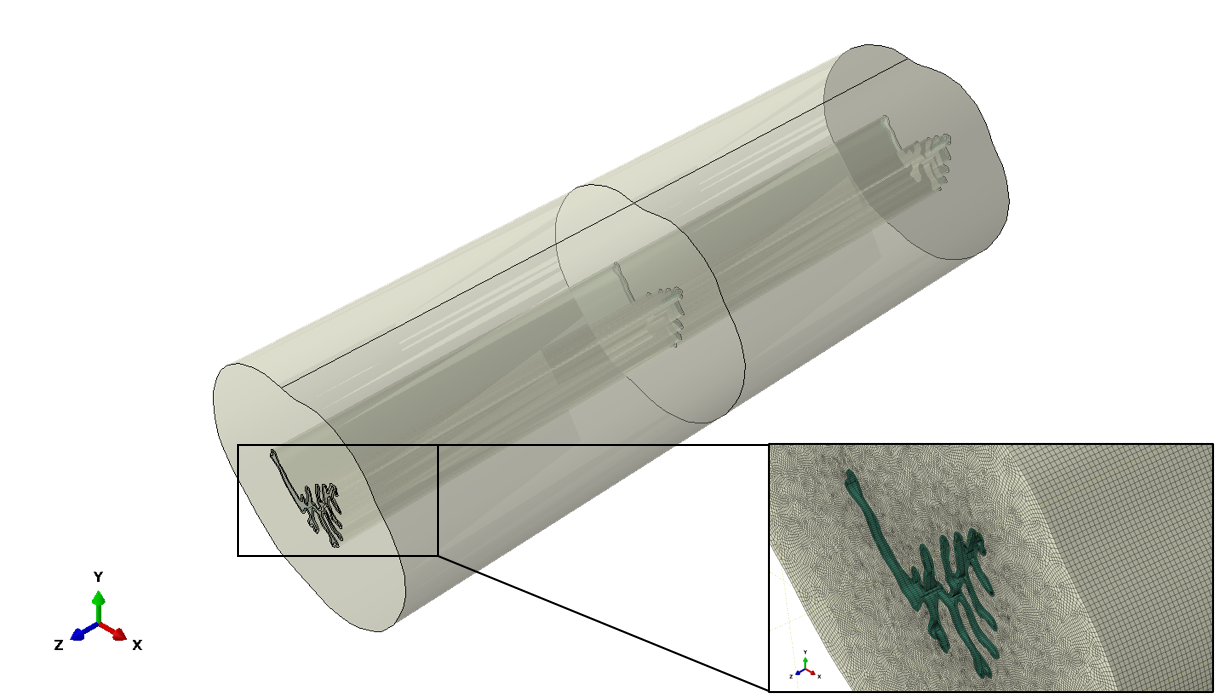}
\caption{Geometry of the test case {\tt c4zz134m}. The box represents the mesh of the
         urethral duct. The two colors refer to the different material of the model.}
\label{fig:c4zz_mesh}
\end{figure}

The mesh derives from the discrtization of the complex conformation of the urethral duct,
with particular regard to the bulbar region. The duct locally consists of an inner thin
layer of dense connective tissue and an outer thick stratum of more compliant spongy
tissue \cite{natali2017urethral,natali2017mechanics}. Both the materials are linear elastic,
characterized by $(E,\nu) = (0.06 MPa, 0.4)$ and $(E,\nu) = (0.0066 MPa, 0.4)$,
respectively. The Figure \ref{fig:c4zz_mesh} shows a representation of the problem's
geometry and mesh.

\bibliographystyle{siamplain}
\bibliography{biblio}

\end{document}